      \documentclass[11pt]{amsart} 

      \usepackage[mathscr]{eucal} 
      \usepackage{amsmath,amsfonts, amssymb,latexsym} 
\def\rg{\hbox to 30pt{\rightarrowfill}}
\def\lg{\hbox to 30pt{\leftarrowfill}}

      \parskip\smallskipamount 
          \newtheorem{theorem}{Theorem}[section] 
       
      \newtheorem{proposition}[theorem]{Proposition} 
      \newtheorem{corollary}[theorem]{Corollary} 
      \newtheorem{lemma}[theorem]{Lemma}

      \makeatletter 
      \@addtoreset{equation}{section} 
      \makeatother 

\newcommand{\HH}{{\mathbb H}}
       
      \newcommand{\CC}{{\mathbb C}}

      \newcommand{\DD}{{\mathbb D}} 
       
      \newcommand{\FF}{{\mathbb F}} 
      \newcommand{\TT}{{\mathbb T}}

      \newcommand{\cA}{{\mathcal A}}

      \newcommand{\cD}{{\mathcal D}}
\newcommand{\cQ}{{\mathcal Q}} 
      \newcommand{\cE}{{\mathcal E}} 
      \newcommand{\cF}{{\mathcal F}} 
      \newcommand{\cG}{{\mathcal G}} 
      \newcommand{\cH}{{\mathcal H}}

\newcommand{\A}{{\mathcal A}}

      \newcommand{\E}{{\mathcal E}} 
      \newcommand{\F}{{\mathcal F}}

             \newcommand{\R}{{\mathcal R}}

      \newcommand{\V}{{\mathcal V}}

      \newcommand{\cK}{{\mathcal K}} 
      \newcommand{\cL}{{\mathcal L}} 
      \newcommand{\cM}{{\mathcal M}} 
      \newcommand{\cN}{{\mathcal N}}

      \newcommand{\cP}{{\mathcal P}} 
      \newcommand{\cR}{{\mathcal R}}

      \newcommand{\cY}{{\mathcal Y}} 
      \newcommand{\cX}{{\mathcal X}} 
      \newcommand{\cW}{{\mathcal W}}

      \newcommand{\rank}{\hbox{\rm{rank}}\,}

      \newdimen\expt 
      \def\boxit#1{\setbox0\hbox{$\displaystyle{#1}$} 
            \hbox{\lower.4\expt 
       \hbox{\lower3\expt\hbox{\lower\dp0 
            \hbox{\vbox{\hrule height.4\expt 
       \hbox{\vrule width.4\expt\hskip3\expt 
            \vbox{\vskip3\expt\box0\vskip2\expt}%
       \hskip3\expt\vrule width.4\expt}\hrule height.4\expt}}}}}} 
      
\begin{document} 
\pagestyle{myheadings}
\markboth{ Gelu Popescu}{ Factorizations of  characteristic functions and joint invariant subspaces }


\title [ Characteristic functions and  joint invariant subspaces ] 
{Characteristic functions and  joint invariant subspaces }
  \author{Gelu Popescu}

\address{Department of Mathematics, The University of Texas 
at San Antonio \\ San Antonio, TX 78249, USA}
\email{\tt gelu.popescu@utsa.edu}

\date{July 1, 2005}
\thanks{Research supported in part by an NSF grant}

\begin{abstract}
Let $T:=[T_1,\ldots, T_n]$ be an $n$-tuple of operators on a  Hilbert space such that $T$ is a completely non-coisometric row contraction.
We establish the existence of a ``one-to-one'' correspondence between the joint invariant subspaces  under $T_1,\ldots, T_n$, and the regular factorizations of the characteristic function $\Theta_T$ associated with $T$. 
In particular, we prove that
  there is a non-trivial joint invariant subspace  under 
the  operators $T_1,\ldots, T_n$, if and only if there is a non-trivial regular factorization of $\Theta_T$. We also provide a functional model for the joint invariant subspaces in terms of the regular  factorizations of the characteristic function, and  prove the existence of joint invariant subspaces
for certain classes
of $n$-tuples of operators. 

We obtain criterions for joint similarity of $n$-tuples of operators to Cuntz row isometries. In particular, we prove that a  completely non-coisometric row contraction $T$ is jointly similar to a Cuntz row isometry if and only if the characteristic function of $T$  is an invertible multi-analytic operator.

\end{abstract}

\maketitle

\section{Introduction}

In the classical case of a single operator, the connection between the invariant subspaces of an operator and the corresponding characteristic function was first considered, for certain particular classes of operators, in the work of Liv\v{s}itz, Potapov, \v{S}mulyan, Brodsky, etc (see the references from
\cite{SzF1} and \cite{SzF2}).
One of the fundamental results in the Nagy-Foia\c s theory of contractions \cite{SzF-book} states that the invariant subspaces of a completely non-unitary (c.n.u.) contraction $T$ on a (separable) Hilbert space are in ``one-to-one'' correspondence with the regular factorizations of the characteristic function  associated with $T$. This general result, although influenced in part by the work  of the authors cited above, was obtained by Sz.-Nagy and Foia\c s in \cite{SzF1},  \cite{SzF2},  following an entirely different approach based on the geometric structure of the unitary dilation  and the corresponding functional model for 
c.n.u. contractions.
 
The main goal of this paper is to  obtain a multivariable version of the above-mentioned result, for $n$-tuples of operators, and to provide a functional model for the joint invariant subspaces in terms of the regular  factorizations of the characteristic function.
This  comes as a natural continuation of our program to develop a {\it free} analogue of Nagy-Foia\c s theory, for row contractions.

An $n$-tuple $T:=[T_1,\ldots, T_n]$ of  bounded linear operators acting on a common Hilbert space $\cH$ is called  {\it row contraction} if
$$
T_1T_1^*+\cdots+T_nT_n^*\leq I.
$$
A distinguished role among row contractions is played by 
the $n$-tuple $S:=[S_1,\ldots, S_n]$ of {\it left creation operators} on the full Fock space with $n$ generators, $F^2(H_n)$, which satisfies the noncommutative von Neumann inequality \cite{Po-von} (see also \cite{Po-funct}, \cite{Po-poisson})
$$\|p(T_1,\ldots, T_n)\|\leq \|p(S_1,\ldots, S_n)\|
$$
for any polynomial $p(X_1,\ldots, X_n)$ in $n$ noncommuting indeterminates. For the classical von Neumann inequality \cite{vN} (case $n=1$) and a nice survey,  we refer to Pisier's book \cite{Pi}.
Based on the left creation operators and their representations, a noncommutative dilation theory and model theory for row contractions was developed in \cite{F}, \cite{B}, \cite{Po-models}, \cite{Po-isometric}, 
\cite{Po-charact},
 \cite{Po-intert}, etc.
In this study, the role of the unilateral shift is played by 
 the left creation operators  and the Hardy algebra $H^\infty(\DD)$
 is replaced by the noncommutative analytic Toeplitz algebra $F_n^\infty$. We recall that $F_n^\infty$  was introduced in \cite{Po-von} as the algebra of left multipliers of  
$F^2(H_n)$  and  can be identified with 
the
 weakly closed  (or $w^*$-closed) algebra generated by the left creation operators
   $S_1,\dots, S_n$      and the identity.

 In \cite{Po-charact}, we defined the {\it standard characteristic function} of a row contraction 
 (a multi-analytic operator acting on Fock spaces) which, as in the classical case ($n=1$)
 \cite{SzF-book}, turned out to be a complete unitary invariant for completely 
 non-coisometric
 row contractions (c.n.c.).  We also constructed a model for c.n.c. row contractions, in which the characteristic function occurs explicitely.

In 2000, Arveson \cite{Arv2} introduced and studied the curvature 
 and Euler characteristic associated with a row contraction with commuting entries. 
 Noncommutative
 analogues of these numerical invariants were defined and studied  by the
 author \cite{Po-curvature} and, independently, by D.~Kribs \cite{Kr}.
  We showed in \cite{Po-varieties} that the curvature invariant and Euler characteristic 
asssociated with a Hilbert module generated by  an arbitrary 
(resp.~commuting) row contraction $T:=[T_1,\ldots, T_n]$ can be expressed only in terms of the 
(resp.~constrained) characteristic function of $T$.

In this paper, we continue the study of the characteristic function $\Theta_T$ associated with a row contraction $T:=[T_1,\ldots, T_n]$ in connection with joint invariant subspaces under the operators $T_1,\ldots, T_n$, and the joint similarity of $T$  to a Cuntz row isometry $W:=[W_1,\ldots, W_n]$, i.e.,  $W_1,\ldots, W_n$ are isometries with 
$$
W_1 W_n^*+\cdots +W_nW_n^*=I.
$$
 After some preliminaries on multivariable  noncommutative dilation theory (see Section 2),  we present in Section 3 the main results of this paper.

  We establish the existence of a ``one-to-one'' correspondence between the joint invariant subspaces  under $T_1,\ldots, T_n$, and the regular factorizations of the characteristic function $\Theta_T$ associated with   a completely non-coisometric row contraction $T:=[T_1,\ldots, T_n]$ (see Theorem \ref{inv-factor1} and Theorem \ref{inv-factor2}). 
In particular, we prove that
  there is a non-trivial joint invariant subspace  under 
the  operators $T_1,\ldots, T_n$, if and only if there is a non-trivial regular factorization of $\Theta_T$ (see Theorem \ref{non-trivial sub}). Using the model theory for c.n.c. row contractions, we   provide a functional model for the joint invariant subspaces in terms of the regular factorizations of the characteristic function (see Theorem \ref{inv-factor}).
An important question related to the main result, Theorem \ref{inv-factor1}, is to what extent a joint invariant subspace determines the corresponding regular factorization
of the characteristic function. We address this problem in Theorem \ref{inv-factor3}.

In Section 4, we prove the existence of a unique triangulation of type
$$
\left(\begin{matrix}C_{\cdot 0}&0\\
*& C_{\cdot 1}\end{matrix}
\right)
$$
for any row contraction $T:=[T_1,\ldots, T_n]$ (see Theorem \ref{factori}), and prove the existence of non-trivial joint invariant subspaces
for certain classes
of  row contractions.
We also  show that there is a non-trivial joint invariant subspace under $T_1,\ldots, T_n$
whenever the inner-outer factorization of the characteristic function  associated with $T $ is non-trivial (see Theorem \ref{chara-fact}).

In Section 5, we obtain criterions for joint similarity of $n$-tuples of operators to Cuntz row isometries. In particular, we prove that a  completely non-coisometric row contraction $T$ is jointly similar to a Cuntz row isometry if and only if the characteristic function of $T$  is an invertible multi-analytic operator (see Theorem \ref{simi-Cuntz}). Moreover, in this case, we provide a model Cuntz row isometry for  similarity.  This is a multivariable version of a result of Sz.-Nagy and Foia\c s \cite{SzF3}, concerning the  similarity to unitary operators.

 Extending on some results
obtained  by Sz.-Nagy   \cite {SzN}, Nagy-Foia\c s  \cite{SzF-book},  and the author \cite{Po-models}, \cite{Po-similarity}, we prove, in particular, that
 a  one-to-one power bounded $n$-tuple $[T_1,\ldots, T_n]$ 
 of operators  on a Hilbert space $\cH$ is jointly similar to a Cuntz row isometry if and only if
there exists  a constant $c>0$
such that
\begin{equation*}  
\sum\limits_{\alpha\in\FF_n^+, |\alpha|=k} \|T_\alpha^*h\|^2\geq c\|h\|^2,\quad h\in\cH, 
\end{equation*}
for any $k=1,2,\ldots$.

Recently \cite{Po-varieties}, \cite{Po-varieties2} we developed a dilation theory for row contractions $[T_1,\ldots, T_n]$ subject to constraints such as ~$p(T_1,\ldots,T_n)=0$,~ ~$p\in \cP$, where $\cP$ is a set of noncommutative polynomials.
It would be interesting to see to what extent the results of this paper can be extended to {\it constrained row contractions} and their {\it constrained characteristic functions}.

\bigskip

\section{Preliminaries on characteristic functions for row contractions}

Let $H_n$ be an $n$-dimensional complex  Hilbert space with orthonormal basis
$e_1$, $e_2$, $\dots,e_n$, where $n\in \{1,2,\dots\}$ or $n=\infty$.
  We consider the full Fock space  of $H_n$ defined by

$$F^2(H_n):=\bigoplus_{k\geq 0} H_n^{\otimes k},$$ 
where $H_n^{\otimes 0}:=\CC 1$ and $H_n^{\otimes k}$ is the (Hilbert)
tensor product of $k$ copies of $H_n$.
Define the left creation 
operators $S_i:F^2(H_n)\to F^2(H_n), \  i=1,\dots, n$,  by
$$
 S_i\varphi:=e_i\otimes\varphi, \quad  \varphi\in F^2(H_n).
$$

 The  noncommutative analytic Toeplitz   algebra   $F_n^\infty$ 
  and  its norm closed version,
  the noncommutative disc
 algebra  $\cA_n$,  were introduced by the author   \cite{Po-von} in connection
   with a multivariable noncommutative von Neumann inequality.
$F_n^\infty$  is the algebra of left multipliers of  
$F^2(H_n)$  and  can be identified with 
 the
  weakly closed  (or $w^*$-closed) algebra generated by the left creation operators
   $S_1,\dots, S_n$  acting on   $F^2(H_n)$,
    and the identity.
     When $n=1$, $F_1^\infty$ 
  can be identified
   with $H^\infty(\DD)$, the algebra of bounded analytic functions
    on the open unit disc. The algebra $F_n^\infty$ can be viewed as a
     multivariable noncommutative 
    analogue of $H^\infty(\DD)$.
 There are many analogies with the invariant
  subspaces of the unilateral 
 shift on $H^2(\DD)$, inner-outer factorizations,
  analytic operators, Toeplitz operators, $H^\infty(\DD)$--functional
   calculus, bounded (resp.~spectral) interpolation, etc.

Let $\FF_n^+$ be the unital free semigroup on $n$ generators 
$g_1,\dots,g_n$, and the identity $g_0$.
The length of $\alpha\in\FF_n^+$ is defined by
$|\alpha|:=k$, if $\alpha=g_{i_1}g_{i_2}\cdots g_{i_k}$, and
$|\alpha|:=0$, if $\alpha=g_0$.
We also define
$e_\alpha :=  e_{i_1}\otimes e_{i_2}\otimes \cdots \otimes e_{i_k}$  
 and $e_{g_0}= 1$.
It is  clear that  $\{e_\alpha:\alpha\in\FF_n^+\}$ is an orthonormal basis of $F^2(H_n)$.
   If $T_1,\dots,T_n\in B(\cH)$, the algebra of all bounded linear operators on a Hilbert space $\cH$, we define 
$T_\alpha :=  T_{i_1}T_{i_2}\cdots T_{i_k}$
 and 
$T_{g_0}:=I_\cH$.

We need to recall from  \cite{Po-charact},
\cite{Po-multi}, \cite{Po-von},  \cite{Po-funct}, and  \cite{Po-analytic} 
 a few facts
 concerning multi-analytic   operators on Fock spaces.
   We say that 
 a bounded linear
  operator 
$A$ acting from $F^2(H_n)\otimes \cK$ to $ F^2(H_n)\otimes \cK'$ is 
 multi-analytic
if 
\begin{equation}
A(S_i\otimes I_\cK)= (S_i\otimes I_{\cK'}) A\quad 
\text{\rm for any }\ i=1,\dots, n.
\end{equation}
Notice that $A$ is uniquely determined by the operator
$\theta:\cK\to F^2(H_n)\otimes \cK'$, which is   defined by ~$\theta k:=A(1\otimes k)$, \ $k\in \cK$, 
 and
is called the  symbol  of  $A$. We denote $A=A_\theta$. Moreover, 
$A_\theta$ is 
 uniquely determined by the ``coefficients'' 
  $\theta_{(\alpha)}\in B(\cK, \cK')$, which are given by
 $$
\left< \theta_{(\tilde\alpha)}k,k'\right>:= \left< \theta k, e_\alpha 
\otimes k'\right>=\left< A_\theta(1\otimes k), e_\alpha 
\otimes k'\right>,\quad 
k\in \cK,\ k'\in \cK',\ \alpha\in \FF_n^+,
$$
where $\tilde\alpha$ is the reverse of $\alpha$, i.e., $\tilde\alpha= g_{i_k}\cdots g_{i_1}$ if
$\alpha= g_{i_1}\cdots g_{i_k}$.
We can associate with $A_\theta$ a unique formal Fourier expansion 

\begin{equation*}
A_\theta\sim \sum_{\alpha \in \FF_n^+} R_\alpha \otimes \theta_{(\alpha)},
\end{equation*}
where $R_i:=U^*S_iU$, \ $i=1,\ldots, n$, are the right creation operators
on $F^2(H_n)$  and $U$ is the unitary operator on $F^2(H_n)$
mapping $e_{i_1}\otimes e_{i_2}\otimes \cdots \otimes e_{i_k}$ into $e_{i_k}\otimes\cdots \otimes e_{i_2}\otimes e_{i_1}$.
Based on the noncommutative von Neumann
   inequality \cite{Po-funct},
  we proved that
  $$A_\theta=\text{\rm SOT}-\lim_{r\to 1}\sum_{k=0}^\infty \sum_{|\alpha|=k}
   r^{|\alpha|} R_\alpha\otimes \theta_{(\alpha)},
   $$
   where, for each $r\in (0,1)$ the series converges in the uniform norm.
The set of  all multi-analytic operators in 
$B(F^2(H_n)\otimes \cK,
F^2(H_n)\otimes \cK')$  coincides  with   
$R_n^\infty\bar \otimes B(\cK,\cK')$, 
the WOT closed algebra generated by the spatial tensor product, where $R_n^\infty:=U^* F_n^\infty U$
(see \cite{Po-analytic}  and \cite{Po-central}).
The multi-analytic operator  $~A_\theta~$ is  called
\begin{enumerate}
\item [(i)]
inner if $~A_\theta~$ is an isometry,
\item [(ii)]
outer if $~\overline{A_\theta (F^2(H_n)\otimes\cE)}=F^2(H_n)\otimes\cE_*~$,
\item [(iii)]
purely contractive if $~\|P_{\E_*}\theta h\|<\|h\|~$ for every $~h\in\cE,~h\ne 0~$,
\item [(iv)] unitary constant if $A_\theta=I\otimes W$ for some unitary operator $W\in B(\cK,\cK')$.
\end{enumerate}
%
If $A_{\theta'}:F^2(H_n)\otimes \cE'\to F^2(H_n)\otimes \cE_*'$
is  another multi-analytic operator, we say that $A_\theta$ coincides
with $A_{\theta'}$ 
if there exist two unitary operators
$$
W:\cE\to\cE'~,\quad W_*:\cE_*\to\cE_*'
$$
such that
$$
(I\otimes W_*)A_\theta=A_{\theta'}(I\otimes W).
$$
For simplicity, throughout this paper,
 $T:=[T_1,\ldots, T_n]$, \ $n=1,\ldots, \infty$,  denotes either the $n$-tuple $(T_1,\ldots, T_n)$  of bounded linear operators on a Hilbert space $\cH$ or the row operator matrix $[T_1 ~\cdots ~T_n]$  acting from $\cH^{(n)} $ to $\cH$, where $\cH^{(n)}:=\oplus_{i=1}^n \cH$  is the direct sum of $n$ copies of $\cH$.
Assume that 
$T:=[T_1,\ldots, T_n]$ is a  row contraction, i.e.,
$$
T_1T_1^*+\cdots +T_nT_n^*\leq I.
$$ 
 The defect operators of $T$ are
$$
\Delta_{T^*}:=\left(I_\cH-\sum_{i=1}^n T_i T_i^*\right)^{1/2}\in B(\cH)\quad \text{ and } \quad \Delta_T:=(I_{\cH^{(n)}}-T^*T)^{1/2}\in B(\cH^{(n)}),
$$
 and the defect spaces of $T$ are defined by
$$
\cD_*:=\overline{\Delta_{T^*}\cH}\quad  \text{ and }\quad \cD:=\overline{\Delta_T \cH^{(n)}}.
$$
The characteristic function of  the row contraction $T:=[T_1,\ldots, T_n]$   is the multi-analytic operator
$\Theta_T: F^2(H_n)\otimes \cD\to F^2(H_n)\otimes \cD_*$
with symbol $\theta_T$ is given by
$$
\theta_T(h):=-\sum_{i=1}^n T_i P_i h+\sum_{i=1}^n (S_i\otimes I_{\cD_*})\left(\sum_{\alpha\in \FF_n^+} e_\alpha \otimes \Delta_{T^*}T_\alpha^*P_i \Delta_T h\right),\qquad h\in\cD,
$$
where $P_i$ denotes the orthogonal projection of $ \cH^{(n)}$ onto the $i$-component of $\cH^{(n)}$,
 and  $S:=[S_1,\ldots, S_n]$  is the model multi-shift  of left creation operators acting on the full Fock space $F^2(H_n)$.

Using the characterization of multi-analytic operators on Fock spaces (see \cite{Po-analytic}, \cite{Po-central}), one can easily see that the characteristic  function  of $T$
is  a multi-analytic operator
with the formal Fourier representation
\begin{equation*}
  -I\otimes T+
\left(I\otimes \Delta_{T^*}\right)\left(I -\sum_{i=1}^n R_i\otimes T_i^*\right)^{-1}
\left[R_1\otimes I_\cH,\ldots, R_n\otimes I_\cH
\right] \left(I\otimes \Delta_{T}\right),
 \end{equation*}
where $R_1,\ldots, R_n$ are the right creation operators on the full Fock space $F^2(H_n)$.

The definition of the characteristic function of $T$ arises in a natural way in the context of the theory of noncommutative isometric dilations  for row contractions (see
\cite{Po-isometric} and \cite{Po-charact}).
Let $V:=[V_1,\ldots, V_n]$, $V_i\in B(\cK)$,  be the minimal isometric dilation of $T$ on a Hilbert space
$\cK\supset\cH$. Therefore, 
\begin{enumerate}
\item[(i)] $V_1,\ldots, V_n$ are isometries with orthogonal ranges;
\item[(ii)] $T_i^*=V_i^*|_\cH$, $i=1,\ldots, n$;
\item[(iii)]
$\cK=\bigvee_{\alpha\in \FF_n^+} V_\alpha \cH$.
\end{enumerate}
Consider the following subspaces of $~\cK$:
$$
\cL:=\bigvee_{i=1}^n(V_i-T_i)\cH,\qquad\cL_*:=\overline{\left(I_\cK-\sum_{i=1}^n V_i T_i^*\right)\cH}.
$$
According to \cite{Po-isometric},
we have the following  orthogonal decompositions of the minimal isometric dilation space of $T$:
\begin{equation}\label{two-dec}
\cK=\R \oplus M_V(\cL_*)=\cH \oplus M_V(\cL),
\end{equation}
where $\cR$ reduces each operator $V_i$, $i=1,\ldots,n$, 
$$
M_V(\cL_*)=\bigoplus_{\alpha\in \FF_n^+} V_\alpha \cL_*,
\quad \text{ and } \quad M_V(\cL)=\bigoplus_{\alpha\in \FF_n^+} V_\alpha \cL.
$$  
Denote by $~\Phi^\cL~$ the unitary operator from $~M_V(\cL)~$ to $~F^2(H_n)\otimes \cL$
defined by
$$
\Phi^\cL\left(\sum\limits_{\alpha\in\FF_n^+}V_\alpha
\ell_\alpha\right):=\sum\limits_{\alpha\in\FF_n^+}e_\alpha \otimes\ell_\alpha,\qquad
\ell_\alpha\in\cL,\ \sum\limits_{\alpha\in\FF_n^+}\|\ell_\alpha\|^2<\infty.
$$
One can view $\Phi^\cL$ as the Fourier representation
of $~M_V(\cL)~$ on Fock spaces.
Then, for any $i=1,\ldots, n$, we have
$$
\Phi^\cL V_i=(S_i\otimes I_\cL)\Phi^\cL,
$$
where  $S:=[S_1,\ldots, S_n]$ is the model multi-shift  of left creation operators acting 
on the full Fock space $F^2(H_n)$. Similarly, one can define 
the unitary operator 
(Fourier representation)
 $\Phi^{\cL_*}: M_V(\cL_*)\to F^2(H_n)\otimes \cL_*$.
 We proved in \cite{Po-charact} that  the characteristic function $\Theta_T$  coincides with  the multi-analytic operator
$~{\Theta_\cL}:F^2(H_n)\otimes \cL\to F^2(H_n)\otimes \cL_*$
defined by
$$
{\Theta_\cL}:=\Phi^{\cL_*}(P_{M_V(\cL_*)}|_{M_V(\cL)})(\Phi^\cL)^*,
$$
where $P_{M_V(\cL_*)}$ denotes the orthogonal projection of $~\cK~$ onto $~M_V(\cL_*)$.
 \bigskip

Let $T:=[T_1,\ldots, T_n]$, \ $n=1,\ldots, \infty$, be a row contraction with $T_i\in B(\cH)$ and consider the subspace $\cH_c\subset \cH$ defined by
$$
\cH_c:=\left\{h\in \cH:\ \sum_{|\alpha|=k} \|T_\alpha^*h\|^2=\|h\|^2
\text{ for any } k=1,2,\ldots\right\}
$$
We call $T$  a completely non-coisometric (c.n.c.)  row contraction if $\cH_c=\{0\}$. We proved in \cite{Po-isometric} that $\cH_c$ is a joint invariant subspaces under the operators $T_1^*,\ldots, T_n^*$, and it is  also the largest subspace in $\cH$ on which $T^*$  acts isometrically. Consequently, we have the following triangulation with respect to the decomposition
$\cH=\cH_c\oplus \cH_{cnc}$:
$$
T_i=\left(\begin{matrix}A_i&0\\
*&B_i
\end{matrix}\right),\qquad i=1,\ldots,n,
$$
where $[A_1,\ldots, A_n]$ is a coisometry, i.e.,
$A_1 A_1^*+\cdots +A_n A_n^*=I_{\cH_c}$, and
$[B_1,\ldots, B_n]$ is a c.n.c. row contraction.

 In \cite{Po-charact}, we  constructed the following model for  c.n.c. row contractions,  in which the characteristic function occurs explicitly. 
\begin{theorem}\label{funct-model1}
Every completely non-coisometric row contraction $T:=[T_1,\ldots, T_n]$,      $n=1,2,\ldots,\infty$,  on a Hilbert space $\cH$ is unitarily equivalent to a row contraction $\bf{T}:=[{\bf T}_1,\ldots, {\bf T}_n]$
 on the Hilbert space
$$
{\bf H}:=[ (F^2(H_n)\otimes \cD_*)\oplus\overline{\Delta_{\Theta_T} (F^2(H_n)\otimes \cD)}]\ominus
\{\Theta_T f\oplus\Delta_{\Theta_T} f:\ \ f\in F^2(H_n)\otimes \cD\},
$$

where $~\Delta_{\Theta_T}:=(I-\Theta_T^*\Theta_T)^{1/2}$ and
  the operator ${\bf T}_i$, $i=1,\ldots, n$, is defined by
$$
{\bf T}_{i}^*[f\oplus\Delta_{\Theta_T} (S_j\otimes I_{\cD_*})  g]:=\begin{cases}
(S_{i}^*\otimes I_{\cD_*}) f\oplus\Delta_{\Theta_T} g&\qquad \text{if } \ i =j, \\
(S_{i}^* \otimes I_{\cD_*})f \oplus 0&\qquad \text{if } \ i \ne j,
\end{cases}
$$
$i,j=1,\ldots,n$, and  $S_1,\ldots, S_n$ are the left creation operators on the full Fock space $F^2(H_n)$.

 Moreover, $~T$ is a pure row contraction if and only  if   $~\Theta_T~$ is an inner multi-analytic operator. In this  case the model reduces to
$$
{\bf H}= (F^2(H_n)\otimes {\cD_*})\ominus \Theta_T (F^2(H_n)\otimes \cD),\qquad{\bf T}_{i}^* f=(S_{i}^*\otimes I_{\cD_*}) f,\qquad
f\in \bf H.
$$
\end{theorem}

 Any contractive multi-analytic operator
$\Theta:F^2(H_n)\otimes \cE\to F^2(H_n)\otimes {\cE_*}\quad(\cE,\cE_*~$ are Hilbert spaces)  generates a c.n.c. row contraction 
${\bf T}:=[{\bf T}_1,\ldots, {\bf T}_n]$. More precisely, we proved in \cite{Po-charact} the following result.
\begin{theorem}\label{funct-model2}
 Let $~\Theta:F^2(H_n)\otimes \cE\to F^2(H_n)\otimes \cE_*$ be a contractive multi-analytic operator and  set
 $~\Delta_\Theta:=(I-\Theta^*\Theta)^{1/2}~$.
Then the row contraction ${\bf T}:=[{\bf T}_1,\ldots, {\bf T}_n]$ defined on the Hilbert space
$$
{\bf H}:=[ (F^2(H_n)\otimes {\E_*})\oplus\overline{\Delta_\Theta (F^2(H_n)\otimes \E)}]\ominus\{\Theta g\oplus\Delta_\Theta g :\quad
g\in F^2(H_n)\otimes \cE\}
$$
by
$$
{\bf T}_i^*(f\oplus \Delta_\Theta g):=(S_i^*\otimes I_{\cE_*}) f \oplus C_i^*(\Delta_\Theta g),\qquad  i=1,\ldots,n,
$$
where each operator $~C_i~$ is defined by 
$$C_i(\Delta_\Theta g):=\Delta_\Theta (S_i\otimes I_\cE) g, \quad g\in F^2(H_n)\otimes \cE,
$$
 and   $S_1,\ldots, S_n$ are the left creation operators on $F^2(H_n)$, is completely non-coisometric.

If $~\Theta~$ is purely contractive and 
\begin{equation*}
\overline{\Delta_\Theta(F^2(H_n)\otimes \cE)}=\overline{\Delta_\Theta((F^2(H_n)\otimes \cE)\ominus \cE)},
\end{equation*}
then $~\Theta~$ coincides with the characteristic function of
the row contraction 
${\bf T}:=[{\bf T}_1,\ldots, {\bf T}_n]$.
In this case, considering $~\bf H~$ as a subspace of
$$
{\bf K}:= (F^2(H_n)\otimes {\E_*})\oplus \overline{\Delta_\Theta (F^2(H_n)\otimes \cE)},
$$
we have that the sequence of operators
${\bf V}:=[{\bf V}_1,\ldots, {\bf V}_n]$
 defined on 
$~\bf K~$ by
$$
{\bf V}_i:=(S_i\otimes I_{\cE_*})\oplus C_i,\qquad i=1,\ldots,n,
$$
is the minimal isometric dilation of 
${\bf T}:=[{\bf T}_1,\ldots, {\bf T}_n]$ .
\end{theorem}

\bigskip

\section{Factorizations of characteristic functions and joint invariant subspaces}

  In this section, we establish the existence of a ``one-to-one'' correspondence between the joint invariant subspaces  under $T_1,\ldots, T_n$, and the regular factorizations of the characteristic function $\Theta_T$ associated with a completely non-coisometric row contraction $T:=[T_1,\ldots, T_n]$. 
In particular, we prove that
  there is a non-trivial joint invariant subspace  under 
the  operators $T_1,\ldots, T_n$, if and only if there is a non-trivial regular factorization of $\Theta_T$. Using the model theory for c.n.c. row contractions, we   provide a functional model for the joint invariant subspaces in terms of the regular factorizations of the characteristic function.

Let $\Theta:F^2(H_n)\otimes\cE\to F^2(H_n)\otimes\cE_*$ be a
contractive  
multi-analytic operator and assume that it has the factorization
$$
\Theta=\Theta_2\Theta_1,
$$
where $\Theta_1:F^2(H_n)\otimes\cE\to F^2(H_n)\otimes\cF$ and 
$\Theta_2:F^2(H_n)\otimes\cF\to F^2(H_n)\otimes\cE_*$
are contractive  multi-analytic operators.
Define  the operator
\begin{equation*}
X_\Theta:\overline{\Delta_\Theta(F^2(H_n)\otimes \cE)}\to 
\overline{\Delta_2(F^2(H_n)\otimes \cF)}\oplus \overline{\Delta_1(F^2(H_n)\otimes \cE)}
\end{equation*}
by setting
\begin{equation} \label{X-reg}
X_\Theta(\Delta_\Theta f):= \Delta_2\Theta_1 f\oplus \Delta_1 f,\qquad f\in F^2(H_n)\otimes \cE,
\end{equation}
where $\Delta_\Theta:=(I-\Theta^* \Theta)^{1/2}$ and
$\Delta_j:=(I-\Theta_j^* \Theta_j)^{1/2}$, $j=1,2$.
  Notice that $X_\Theta$ is an isometry. Indeed, 
since 
\begin{equation*}
\begin{split}
I-\Theta^* \Theta&= I-\Theta_1^* \Theta_2^* \Theta_2 \Theta_1\\
&=
\Theta^*_1(I-\Theta_2^* \Theta_2)\Theta_1+(I-\Theta_1^* \Theta_1),
\end{split}
\end{equation*}
we have
\begin{equation*}
\begin{split}
\left\|\Delta_2 \Theta_1 f\oplus \Delta_1 f\right\|^2&=
\|\Delta_2 \Theta_1f\|^2+\|\Delta_1f\|^2\\
&=\left<\Theta^*_1(I-\Theta_2^* \Theta_2)\Theta_1 f+I-\Theta_1^* \Theta_1 f,f\right>\\
&\left(I-\Theta^* \Theta)f,f\right>=\|\Delta_\Theta f\|^2.
\end{split}
\end{equation*}
As in the classical case (see \cite{SzF-book}), we say that the factorization 
$
\Theta=\Theta_2\Theta_1
$
is {\it regular} if $X_\Theta$ is a unitary operator, i.e.,
$$
\left\{ \Delta_2 \Theta_1 f\oplus \Delta_1 f:\ f\in F^2(H_n)\otimes \cE\right\}^{-}=
\overline{\Delta_2(F^2(H_n)\otimes \cF)}\oplus \overline{\Delta_1(F^2(H_n)\otimes \cE)}.
$$

Now let us prove the following technical result which will be very useful in what follows.

\begin{lemma}\label{intert-AB} Let $\Theta:F^2(H_n)\otimes\cE\to F^2(H_n)\otimes\cE_*$ be a contractive
multi-analytic operator and let $C:=[C_1,\ldots, C_n]$ be the row isometry defined  on 
$\overline{\Delta_\Theta(F^2(H_n)\otimes \cE)}$ by setting
$$
C_i\Delta_\Theta f:= \Delta_\Theta(S_i\otimes I_\cE)f,\quad f\in F^2(H_n)\otimes \cE,
$$
for each $i=1,\ldots, n$, where $\Delta_\Theta:=(I-\Theta^* \Theta)^{1/2}$. Then $C$ is a Cuntz row isometry, i.e., 
$C_1C_1^*+\cdots +C_nC_n^*=I$, if and only if 
\begin{equation}\label{over-delta}
\overline{\Delta_\Theta(F^2(H_n)\otimes \cE)}=\overline{\Delta_\Theta((F^2(H_n)\otimes \cE)\ominus \cE)}.
\end{equation}

Assume that $\Theta$ has the factorization
$$
\Theta=\Theta_2\Theta_1,
$$
where $\Theta_1:F^2(H_n)\otimes\cE\to F^2(H_n)\otimes\cF$ and 
$\Theta_2:F^2(H_n)\otimes\cF\to F^2(H_n)\otimes\cE_*$
are contractive multi-analytic operators and let $E:=[E_1,\ldots, E_n]$ and $F:=[F_1,\ldots, F_n]$  be  the corresponding  row isometries  defined on $\overline{\Delta_1(F^2(H_n)\otimes \cE)}$ and 
$\overline{\Delta_2(F^2(H_n)\otimes \cF)}$, respectively.
Then
\begin{equation} \label{Int-XT}
X_\Theta C_i=\left(\begin{matrix} F_i&0\\
0&E_i\end{matrix}\right)X_\Theta,\quad i=1,\ldots,n,
\end{equation}
where the operator $X_\Theta$ is defined by relation \eqref{X-reg}.
Moreover, if the factorization $\Theta=\Theta_2\Theta_1$ is regular, then  $C$ is a Cuntz  row isometry if and only if $E$ and $F$ are Cuntz  row isometries.  
\end{lemma}

\begin{proof}
First, notice that since $\Theta$ is a multi-analytic operator, i.e.,
$$\Theta(S_i\otimes I_\cE)=(S_i\otimes I_{\cE_*})\Theta,\quad i=1,\ldots, n,
$$
we have
\begin{equation*}
\begin{split}
\left< C_i\Delta_\Theta f, C_j\Delta_\Theta g\right>&=
\left< (S_j^*\otimes I_\cE)(I-\Theta^*\Theta) (S_i\otimes I_\cE)f, g\right>\\
&=\left< \delta_{ij}(I-\Theta^* \Theta)f,g\right>=
\delta_{ij}\left<\Delta_\Theta f, \Delta_\Theta g\right>
\end{split}
\end{equation*}
for any $f,g\in F^2(H_n)\otimes \cE$ and $i,j=1,\ldots,n$.
This shows that  the operators $C_1,\ldots, C_n$ are isometries with orthogonal spaces.  Due to the definition of $C_i$, it is clear that $C_1C_1^*+\cdots +C_nC_n^*=I$ if and only if the range of the operator
$[C_1,\ldots, C_n]$ coincides with $\overline{\Delta_\Theta(F^2(H_n)\otimes \cE)}$, which is equivalent
to \eqref{over-delta}.

On the other hand, for each $i=1,\ldots, n$, and $f\in F^2(H_n)\otimes E$, we have
\begin{equation*}
\begin{split}
X_\Theta C_i(\Delta_\Theta f)&=
X_\Theta  \Delta_\Theta(S_i\otimes I_\cE)f \\
&=
\Delta_2\Theta_1(S_i\otimes I_\cE) f\oplus \Delta_1(S_i\otimes I_cE) f\\
&=\Delta_2(S_i\otimes I_\cF)\Theta_1 f\oplus \Delta_1(S_i\otimes I_\cE) f\\
&=F_i\Delta_2\Theta_1 f\oplus E_i\Delta_1 f\\
&=\left(\begin{matrix}
F_i&0\\
0& E_i\end{matrix}\right)\left(\Delta_2 \Theta_1 f \oplus \Delta_1 f\right)\\
&=
\left(\begin{matrix}
F_i&0\\
0& E_i\end{matrix}\right)X_\Theta \Delta_\Theta f,
\end{split}
\end{equation*}
which proves  relation \eqref{Int-XT}.
If  the factorization $\Theta=\Theta_2\Theta_1$ is regular, then $X_\Theta$ is a unitary operator. Consequently, we have
$$
X_\Theta \left(\sum_{i=1}^n C_i C_i^*\right) X_\Theta^*=
\left(\begin{matrix}
\sum\limits_{i=1}^n F_i F_i^* &0\\
0&
\sum\limits_{i=1}^n E_i E_i^*
\end{matrix}
\right),
$$
which implies that $C:=[C_1,\ldots, C_n]$ is a Cuntz  row isometry if and only if $E:=[E_1,\ldots, E_n]$ and $F:=[F_1,\ldots, F_n]$ are Cuntz  row isometries.
This completes the proof.
 \end{proof}

 The main result of this section is the following.

\begin{theorem} \label{inv-factor1}
Let $T:=[T_1,\ldots, T_n]$, $T_i\in B(\cH)$, be a completely non-coisometric row
contraction and let $\Theta:F^2(H_n)\otimes \cE\to F^2(H_n)\otimes \cE_*$ be a contractive multi-analytic operator which coincides with the characteristic function of $T$.
If $\cH_1\subset \cH$ is a joint invariant subspace under the operators $T_1,\ldots, T_n$, then there exists a regular factorization $\Theta=\Theta_2 \Theta_1$, where $\Theta_1:F^2(H_n)\otimes \cE\to F^2(H_n)\otimes \cF$ and 
$\Theta_2:F^2(H_n)\otimes \cF\to F^2(H_n)\otimes \cE_*$
are contractive multi-analytic operators such that
$T:=[T_1,\ldots, T_n]$ is unitarily equivalent to a row contraction
${\bf \TT}:=[{\bf \TT}_1,\ldots, {\bf \TT}_n]$ defined on the Hilbert space
\begin{equation*}
\begin{split}
{\bf \HH}:=[ (F^2(H_n)\otimes {\E_*})&\oplus\overline{\Delta_2 (F^2(H_n)\otimes \F)}\oplus \overline{\Delta_1 (F^2(H_n)\otimes \E)}]\\
&\ominus\{\Theta_2\Theta_1 f\oplus\Delta_2\Theta_1f\oplus \Delta_1f  : \ 
f\in F^2(H_n)\otimes \cE\},
\end{split}
\end{equation*}
   by setting
$$
{\bf \TT}_i^*(f\oplus \varphi\oplus \psi):=(S_i^*\otimes I_{\cE_*}) f \oplus F_i^*\varphi\oplus E_i^*\psi,\qquad  
f\oplus \varphi\oplus \psi\in \HH,
$$
for any $i=1,\ldots,n$, where the operators $F_i$ and $E_i$ are defined in Lemma $\ref{intert-AB}$ and $S_1,\ldots, S_n$ are the left creation operators on $F^2(H_n)$.
Moreover, the subspaces corresponding to $\cH_1$ and $\cH_2:=\cH\ominus \cH_1$ are
\begin{equation*}
\begin{split}
{\bf \HH}_1:=\{\Theta_2 f\oplus \Delta_2f\oplus g&:\ f\in F^2(H_n)\otimes \cF, \ g\in 
\overline{\Delta_1 (F^2(H_n)\otimes \cE)}\}\\
&\ominus
\left\{ \Theta_2 \Theta_1 f\oplus \Delta_2 \Theta_1f \oplus \Delta_1 f:\ f\in F^2(H_n)\otimes \cE\right\}
\end{split}
\end{equation*}
and
\begin{equation*}
\begin{split}
{\bf \HH}_2:=[(F^2(H_n)\otimes \cE_*)&\oplus 
\overline{\Delta_2 (F^2(H_n)\otimes \cF)}\oplus \{0\}]\\
&\ominus\left\{ \Theta_2 f\oplus \Delta_2 f\oplus \{0\}):\ f\in F^2(H_n)\otimes \cF\right\},
\end{split}
\end{equation*}
respectively.
Conversely, every regular factorization $\Theta=\Theta_2 \Theta_1$ generates via the above formulas the subspaces $\HH_1$ and $\HH_2$ with the following properties:
\begin{enumerate}
\item[(i)] $\HH_1$ is invariant under each operator  $\TT_i$, $i=1,\ldots,n$;
\item[(ii)]
$\HH_2=\HH\ominus \HH_1$.
\end{enumerate}
Under the above identification, $\HH_1$ corresponds to a subspace $\cH_1\subset \cH$ which is invariant under each operator $T_i$, $i=1,\ldots,n$.
  \end{theorem}

\begin{proof}
{\bf Part I.}  ~
Let $T:=[T_1,\ldots, T_n]$,\ $T_i\in B(\cH)$, be a row contraction
and let 
$V:=[V_1,\ldots, V_n]$, \ $V_i\in B(\cK)$, be  its the minimal isometric dilation  on a Hilbert  space $\cK=\bigvee_{\alpha\in \FF_n^+} V_\alpha \cH$. Since $V_1,\ldots, V_n$ are isometries with orthogonal ranges,  the noncommutative Wold decomposition \cite{Po-isometric}
provides the orthogonal decomposition
\begin{equation}\label{wold1}
\cK=\cR \oplus M_V(\cL_*),
\end{equation}
where
\begin{equation*}
\cR:=\bigcap_{k=0}^\infty \left[\bigoplus_{|\alpha|=k} V_\alpha \cK\right] \quad \text{ and }\quad
\cL_*:=\overline{\left(I_\cK-\sum_{i=1}^n V_iT_i^*\right)\cH}.
\end{equation*}
Moreover, 
$\cR$ is the maximal subspace of $\cK$   which is reducing for the operators
$V_1,\ldots, V_n$ and the row contraction $[V_1|_{\cR},\ldots, V_n|_{\cR}]$
is a Cuntz row isometry.

Let $\cH_1\subset \cH$ be an invariant subspace under the operators
$T_1, \ldots,  T_n$.
Since $V_i^*|_\cH=T_i^*$, $i=1,\ldots, n$, we deduce that the subspace $\cH_2:=\cH\ominus \cH_1$ is invariant under the operators
$V_1^*,\ldots, V_n^*$.
Therefore, the subspace $\cG:=\cK\ominus \cH_2$ is invariant under
 $V_1,\ldots, V_n$.
Applying  again the noncommutative Wold decomposition  to the row 
isometry
$[V_1|_\cG,\ldots, V_n|_\cG]$, we obtain the orthogonal decomposition
\begin{equation}\label{G-wold}
\cG=\cR_1\oplus M_V(\cQ),
\end{equation}
where
\begin{equation*}
\cR_1:=\bigcap_{k=0}^\infty \left[\bigoplus_{|\alpha|=k} V_\alpha \cG\right] \quad
\text{ and } \quad 
\cQ:=\cG\ominus \left(\bigoplus_{i=1}^n V_i \cG\right).
\end{equation*}
Since $\cR_1$ reduces the operators $V_1,\ldots, V_n$ and $[V_1|_{\cR_1},\ldots, V_n|_{\cR_1}]$ is a Cuntz row isometry, we deduce that
$\cR_1\subset \cR$.
Notice that $\cR_2:=\cR\ominus \cR_1$  is also  a reducing subspace for  $V_1,\ldots, V_n$ and 
$[V_1|_{\cR_2},\ldots, V_n|_{\cR_2}]$ is a Cuntz row isometry.
Using relations \eqref{wold1} and \eqref{G-wold}, we infer that
\begin{equation*}\begin{split}
\cH_2&=\cK\ominus \cG\\
&=
\left[\cR\oplus M_V(\cL_*)\right]\ominus \left[ \cR_1\oplus M_V(\cQ)\right]\\
&=\left[\cR_2\oplus M_V(\cL_*)\right]\ominus M_V(\cQ).
\end{split}
\end{equation*}
Hence, we deduce that
\begin{equation}
\label{inc1}
M_V(\cQ)\subset  \cR_2\oplus M_V(\cL_*).
\end{equation}
On the other hand, due to \eqref{two-dec}, we  have
$$
\cK=\cR\oplus M_V(\cL_*)=\cH\oplus M_V(\cL).
$$
Hence, we obtain
$$\cH=\left[ \cR\oplus M_V(\cL_*)\right]\ominus M_V(\cL).
$$
Since $\cH_2\subset \cH$, the above representations of $\cH$ and $\cH_2$ imply
$$
\left[ \cR_2\oplus M_V(\cL_*)\right] \ominus M_V(\cQ)\subset 
\left[ \cR\oplus M_V(\cL_*)\right]\ominus M_V(\cL).
$$
Taking into account that $\cR=\cR_1\oplus \cR_2$, we have
$$
\left[ \cR_2\oplus M_V(\cL_*)\right] \ominus M_V(\cQ)=
\left[ \cR\oplus M_V(\cL_*)\right]\ominus \left[ \cR_1\oplus M_V(\cQ)\right].
$$
Consequently, we deduce that
\begin{equation}
\label{inc2}
M_V(\cL)\subset \cR_1\oplus M_V(\cQ)
\end{equation}
and
\begin{equation*}\begin{split}
\cH_1&=\cH\ominus \cH_2\\
&=\left[ \cR_1\oplus M_V(\cQ)\right]\ominus M_V(\cL)\\
&=\cG\ominus M_V(\cL).
\end{split}
\end{equation*}

Let $P_{M_V(\cL_*)}$, $P_{M_V(\cQ)}$, $P_\cR$, $P_{\cR_1}$, and 
$P_{\cR_2}$ be the orthogonal projections onto the corresponding spaces.
According to  relations \eqref{inc1} and \eqref{inc2}, for any $x\in M_V(\cQ)$ and $y\in M_V(\cL)$, we have
\begin{equation}
\label{proj}
x= P_{\cR_2} x+P_{M_V(\cL_*)}x \quad \text{ and }\quad y=P_{\cR_1} y +P_{M_V(\cQ)}y.
\end{equation}
In particular, if $x:=P_{M_V(\cQ)}y$ and $y\in M_V(\cL)$, we deduce that
\begin{equation}
\label{y-proj}
y=P_{\cR_1} y +P_{\cR_2}P_{M_V(\cQ)}y+ P_{M_V(\cL_*)}P_{M_V(\cQ)}y.
\end{equation}
Hence and taking into account that the subspace $\cR_1\oplus \cR_2=\cR$ is orthogonal to $M_V(\cL_*)$, we deduce that
\begin{equation}
\label{two eq}
P_{M_V(\cL_*)}y=P_{M_V(\cL_*)}P_{M_V(\cQ)}y\quad  \text{ and }\quad
P_\cR y= P_{\cR_1} y +P_{\cR_2}P_{M_V(\cQ)}y 
\end{equation}
for any $y\in M_V(\cL)$.
Due to relation \eqref{wold1}, we have
\begin{equation}
\label{pr}
P_{\cR}f=\left(I-P_{M_V(\cL_*)}\right)f,\quad f\in \cK.
\end{equation}
On the other hand, relations \eqref{inc2} and \eqref{inc1}
imply
\begin{equation}
\label{pr1}
P_{\cR_1}y=\left(I-P_{M_V(\cQ)}\right)y,\quad y\in M_V(\cL)
\end{equation}
and
\begin{equation}
\label{pr2}
P_{\cR_2}x=\left(I-P_{M_V(\cL_*)}\right)x,\quad x\in M_V(\cQ).
\end{equation}

Assume now that $[T_1,\ldots, T_n]$ is a c.n.c. row contraction. In this case, we  have (see \cite{Po-isometric})
$$
\cK=M_V(\cL)\bigvee M_V(\cL_*)=\cR\oplus M_V(\cL_*), 
$$
which implies
\begin{equation}\label{Pr}
\overline{P_\cR M_V(\cL)}=
\overline{\left(I-P_{M_V(\cL_*)}\right)M_V(\cL)}=\cR.
\end{equation}
Hence and using the second relation in \eqref{two eq}, we deduce that
$$
\overline{P_{\cR_1} M_V(\cL)}=\cR_1\quad \text{ and } \quad
\overline{P_{\cR_2}  P_{M_V(\cQ)}M_V(\cL)}=\cR_2,
$$
and, consequently,
\begin{equation}\label{densities}
\overline{P_{\cR_1} M_V(\cL)}=\cR_1\quad \text{ and } \quad
\overline{P_{\cR_2}  M_V(\cQ)}=\cR_2.
\end{equation}

{\bf Part II.}
Consider the following  contractions:
\begin{equation*}\begin{split}
Q&:=P_{M_V(\cL_*)}|_{M_V(\cL)}: M_V(\cL)\to M_V(\cL_*),\\
Q_1&:=P_{M_V(\cQ)}|_{M_V(\cL)}: M_V(\cL)\to M_V(\cQ),\text{ and}\\
Q_2&:=P_{M_V(\cL_*)}|_{M_V(\cQ)}: M_V(\cQ)\to M_V(\cL_*).
\end{split}
\end{equation*}
Since $M_V(\cL_*)$, $M_V(\cL)$, and $M_V(\cQ)$ are reducing subspaces for the operators $V_1,\ldots, V_n$, we deduce that,
for each $i=1,\ldots, n$, 
\begin{equation*}\begin{split}
Q\left(V_i|_{M_V(\cL)}\right)&=
\left(V_i|_{M_V(\cL_*)}\right)Q,\\
Q_1\left(V_i|_{M_V(\cL)}\right)&=
\left(V_i|_{M_V(\cQ)}\right)Q_1,\text{ and}\\
Q_2\left(V_i|_{M_V(\cQ)}\right)&=
\left(V_i|_{M_V(\cL_*)}\right)Q_2.
\end{split}
\end{equation*}

Let $\Phi^{\cL_*}:M_V(\cL_*)\to F^2(H_n)\otimes \cL_*$ be the Fourier
representation of the subspace $M_V(\cL_*)$, i.e.,
$$
\Phi^{\cL_*}\left( \sum_{\alpha\in \FF_n^+} V_\alpha \ell_\alpha\right):=\sum_{\alpha\in \FF_n^+} e_\alpha \otimes\ell_\alpha,
$$
where $\ell_\alpha\in \cL_*$ and $\sum\limits_{\alpha\in \FF_n^+} \|\ell_\alpha\|^2<\infty$.
Notice that 
$$
\Phi^{\cL_*}\left(V_i|_{M_V(\cL_*)}\right)=\left( S_i\otimes I_{\cL_*}\right)\Phi^{\cL_*},\quad i=1,\ldots,n,
$$
where $S_1,\ldots, S_n$ are the left creation operators on $F^2(H_n)$.
Similarly, we define the Fourier representations of the subspaces 
$M_V(\cL)$ and $M_V(\cQ)$, respectively.
Now,  due to the above intertwining relations satisfied by $Q$, $Q_1$, and $Q_2$, the operators
\begin{equation}\label{Theta12}
\begin{split}
\Theta_\cL&:F^2(H_n)\otimes \cL\to F^2(H_n)\otimes \cL_*,\quad \Theta_\cL:=\Phi^{\cL_*} Q (\Phi^\cL)^*,\\
\Psi_1&:F^2(H_n)\otimes \cL\to F^2(H_n)\otimes \cQ,
\quad \Psi_1:=\Phi^{\cQ} Q_1 (\Phi^\cL)^*,\ 
\text{ and}\\
\Psi_2&:F^2(H_n)\otimes \cQ\to F^2(H_n)\otimes \cL_*,
\quad \Psi_2:=\Phi^{\cL_*} Q_2 (\Phi^\cQ)^*
\end{split}
\end{equation}
are contractive and multi-analytic.
Hence and using  the first equation in \eqref{two eq}, we have

\begin{equation*}
\begin{split}
\Theta_\cL&=\Phi^{\cL_*} Q (\Phi^\cL)^*
=   \Phi^{\cL_*} \left(P_{M_V(\cL_*)}|_{M_V(\cL)}
\right) (\Phi^\cL)^*\\
&=
\Phi^{\cL_*} \left(P_{M_V(\cL_*)}P_{M_V(\cQ)}|_{M_V(\cL)}
\right) (\Phi^\cL)^*\\
&=
\left[\Phi^{\cL_*} \left(P_{M_V(\cL_*)}|_{M_V(\cQ)}\right) (\Phi^\cQ)^*\right]
\left[\Phi^\cQ\left(
P_{M_V(\cQ)}|_{M_V(\cL)}
\right) (\Phi^\cL)^*\right]\\
&=
\left[\Phi^{\cL_*}  Q_2 (\Phi^\cQ)^*\right]
\left[\Phi^\cQ Q_1 (\Phi^\cL)^*\right]\\
&=\Psi_2 \Psi_1.
\end{split}
\end{equation*}

Due to \eqref{pr} and \eqref{Pr}, there exists a unique
unitary operator
$\Phi_\cR:\cR\to \overline{ \Delta_\cL(F^2(H_n)\otimes \cL)}$
such that
\begin{equation}\label{Phi_R}
\Phi_\cR P_\cR\psi:=\Delta_\cL \Phi^\cL \psi,\quad \psi\in M_V(\cL),
\end{equation}
where $\Delta_\cL:=\left(I-\Theta_\cL^* \Theta_\cL\right)^{1/2}$.
Indeed, we have
\begin{equation*}\begin{split}
\|(I-P_{M_V(\cL_*)})\psi\|^2&=\|\psi\|^2-\|P_{M_V(\cL_*)}\psi\|^2\\
&\|\Phi^\cL\psi\|^2-\|\Phi^{\cL_*}P_{M_V(\cL_*)}\psi\|^2\\
&=\|\Phi^\cL\psi\|^2-\|\Theta_\cL\Phi^\cL \psi\|^2\\
&=\|\Delta_\cL\Phi^\cL\psi\|^2.
\end{split}
\end{equation*}
Consequently,
\begin{equation}\label{PPP}
\Phi:=\Phi^{\cL_*}\oplus\Phi_\cR
\end{equation}
is a unitary operator from the  dilation space $~\cK=M_V(\cL_*)\oplus\R~$ onto the Hilbert space
$$
\widetilde{\bf K}:= (F^2(H_n)\otimes {\cL_*})\oplus \overline{\Delta_\cL (F^2(H_n)\otimes \cL)}.
$$
The image of the space $\cH=\cK\ominus M_V(\cL)~$ under the operator $~\Phi$ 
is 
$$
\Phi\cH=\widetilde{\bf H}:=[ (F^2(H_n)\otimes {\cL_*})\oplus\overline{\Delta_\cL (F^2(H_n)\otimes \cL})]\ominus
\{\Theta_\cL f\oplus \Delta_\cL f\ :\ \ f\in F^2(H_n)\otimes \cL\}.
$$
The row contraction ${T}:=[{T}_1,\ldots, {T}_n]$ is transformed under the unitary operator
$\Phi$ into the row contraction
 $\widetilde{\bf T}:=[\widetilde{\bf T}_1,\ldots, \widetilde{\bf T}_n]$, where
$$
\widetilde{\bf T}_i^*(f\oplus \Delta_\cL g):=(S_i^*\otimes I_{\cL_*}) f \oplus \widetilde{C}_i^*(\Delta_\cL g),\qquad  i=1,\ldots,n,
$$
and  each operator $~\widetilde{C}_i~$ is defined by 
$$\widetilde{C}_i(\Delta_\cL g)=\Delta_\cL (S_i\otimes I_\cL) g, \quad g\in F^2(H_n)\otimes \cL.
$$

Notice that, using relations \eqref{pr1}, \eqref{pr2}, and \eqref{densities}, one can show that there are some unitary operators
\begin{equation*}
\Phi_{\cR_1}:\cR_1\to \overline{ \Delta_{\Psi_1}(F^2(H_n)\otimes \cL)}\quad \text{ and } \quad 
\Phi_{\cR_2}:\cR_2\to \overline{ \Delta_{\Psi_2}(F^2(H_n)\otimes \cQ)}
\end{equation*}
uniquely defined by the relations
\begin{equation}\label{pr12}
\begin{split}\Phi_{\cR_1} P_{\cR_1}x:&=\Delta_{\Psi_1} \Phi^\cL x,\quad x\in M_V(\cL),\\
\Phi_{\cR_2} P_{\cR_2}y:&=\Delta_{\Psi_2} \Phi^\cQ y,\quad y\in M_V(\cQ),
\end{split}
\end{equation}  
where $\Delta_{\Psi_j}:=\left(I-\Psi_j^* \Psi_j\right)^{1/2}$ for  $j=1,2$.
Consequently, since $\cR=\cR_2\oplus \cR_1$ and due to relation \eqref{Phi_R}, the operator
$$
X_\cL: \overline{ \Delta_\cL(F^2(H_n)\otimes \cL)}\to
\overline{ \Delta_{\Psi_2}(F^2(H_n)\otimes \cQ)}
\oplus
\overline{ \Delta_{\Psi_1}(F^2(H_n)\otimes \cL)} 
$$
defined by 
\begin{equation}
\label{uni-x}
X_\cL:=\left(\Phi_{\cR_2}\oplus \Phi_{\cR_1}\right) \Phi_\cR^*
\end{equation}
is unitary.
Due to relations \eqref{Phi_R}, \eqref{two eq}, \eqref{pr12}, and \eqref{Theta12},
we deduce that
\begin{equation*}
\begin{split}
X_\cL\Delta_\cL \Phi^\cL y&= X_\cL\Phi_\cR P_\cR y=\left(\Phi_{\cR_2}\oplus\Phi_{\cR_1}\right) P_\cR y\\
&=\left(\Phi_{\cR_2}\oplus\Phi_{\cR_1}\right)\left(P_{\cR_2} P_{M_V(\cQ)}y\oplus P_{\cR_1}y\right)\\
&=\Delta_{\Psi_2}\Phi^\cQ P_{M_V(\cQ)}y\oplus \Delta_{\Psi_1} \Phi^\cL y\\
&=\Delta_{\Psi_2} \Psi_1\Phi^\cL y\oplus \Delta_{\Psi_1} \Phi^\cL y
\end{split}
\end{equation*}
 for any $y\in M_V(\cL)$.
Hence, we have
\begin{equation}
\label{XDE}
X_\cL\Delta_\cL f=\Delta_{\Psi_2} \Psi_1 f \oplus \Delta_{\Psi_1} f,\quad f\in F^2(H_n)\otimes \cL.
\end{equation}
Since $X_\cL$ is  a unitary operator, we also deduce that

$$\left\{ \Delta_{\Psi_2} \Theta_1 f \oplus \Delta_{\Psi_1} f,\quad f\in F^2(H_n)\otimes \cL\right\}^-=\overline{ \Delta_{\Psi_2}(F^2(H_n)\otimes \cQ)}
\bigoplus \overline{ \Delta_{\Psi_1}(F^2(H_n)\otimes \cL)}.
$$
 Due to \eqref{PPP} and \eqref{uni-x}, we have
$$
\Phi=\Phi^{\cL_*}\oplus X_\cL^*\left(\Phi_{\cR_2}\oplus\Phi_{\cR_1}\right).
$$
Now, we need to find the images $\widetilde{\bf H}_1$ and $\widetilde{\bf H}_2$ of $\cH_1$ and $\cH_2$, respectively,  under the unitary operator $\Phi$.
To find $\widetilde{\bf H}_2$, notice first that, due to  relation \eqref{uni-x}, we have
\begin{equation}
\label{Rz}
\Phi_\cR z=X_\cL^*\left(\Phi_{\cR_2}\oplus\Phi_{\cR_1}\right)(z\oplus 0)=X_\cL^*\left( \Phi_{\cR_2}z\oplus 0\right)
\end{equation}
for any $z\in \cR_2$.
Hence and using \eqref{Phi_R}, we infer that
\begin{equation*}
\begin{split}
\Phi\left( M_V(\cL_*)\oplus \cR_2 \right)&=
\Phi^{\cL_*} M_V(\cL_*)\oplus \Phi_\cR \cR_2\\
&=
(F^2(H_n)\otimes \cL_*)\bigoplus X_\cL^*\left( \overline{ \Delta_{\Psi_2}(F^2(H_n)\otimes \cQ)}\oplus \{0\}\right)
\end{split}
\end{equation*}
and, due to \eqref{proj},
$$
\Phi M_V(\cQ)=\left\{\Phi^{\cL_*} P_{M_V(\cL_*)}f\oplus \Phi_\cR P_{\cR_2} f:\  f\in M_V(\cQ)\right\}.
$$
Hence, and using relations \eqref{Theta12}, \eqref{pr12}, 
and \eqref{Rz}, we obtain
$$
\Phi M_V(\cQ)=\left\{ \Psi_2 u\oplus X_\cL^*(\Delta_{\Psi_2} u\oplus 0):\ u\in F^2(H_n)\otimes \cQ\right\}.
$$
Now, using the representation of $\cH_2$ from Part I, i.e.,
$$
{\cH}_2=\left[M_V(\cL_*)\oplus \cR_2\right]\ominus M_V(\cQ),
$$
 we obtain 
 \begin{equation*}
\begin{split} \widetilde{\bf H}_2=[(F^2(H_n)\otimes \cL_*)&\bigoplus X_\cL^{*}
\left(
\overline{\Delta_{\Psi_2} (F^2(H_n)\otimes \cQ)})\oplus \{0\}\right)]\\
&\ominus\left\{ \Psi_2 f\oplus X_\cL^{*}(\Delta_{\Psi_2} f\oplus 0):\ f\in F^2(H_n)\otimes \cQ\right\}.
\end{split}
\end{equation*}
Since $\widetilde{\bf H}_1=\widetilde{\bf H}\ominus \widetilde{\bf H}_2$, we deduce that
\begin{equation*} 
\begin{split}
\widetilde{\bf H}_1=\{\Psi_2 f\oplus X_\cL^{*}(\Delta_{\Psi_2}f\oplus g)&:\ f\in F^2(H_n)\otimes \cQ,\  g\in 
\overline{\Delta_{\Psi_1}  (F^2(H_n)\otimes \cL)}\}\\
&\ominus
\left\{ \Theta_\cL w\oplus \Delta_\Theta w:\ w\in F^2(H_n)\otimes \cL\right\}.
\end{split}
\end{equation*}

According to  Section 2, the characteristic function  $\Theta_T$ of the row contraction $T$ coincides with $\Theta_\cL$, and therefore with $\Theta$. Via this identification, the regular factorization 
$\Theta_\cL=\Psi_2 \Psi_1$ corresponds to a  regular factorization 
$\Theta=\Theta_2 \Theta_1$,
where  $\Theta_1:F^2(H_n)\otimes\cE\to F^2(H_n)\otimes\cF$ and 
$\Theta_2:F^2(H_n)\otimes\cF\to F^2(H_n)\otimes\cE_*$
are contractive  multi-analytic operators.
 Now, it is easy to see that, under the above identification,  the subspaces $\widetilde{\bf H}_1$ and $\widetilde{\bf H}_2$ correspond to the subspaces
  \begin{equation}
\begin{split}\label{H2}
{\bf H}_2=[(F^2(H_n)\otimes \cE_*)&\oplus X_\Theta^{*}
\left(
\overline{\Delta_2 (F^2(H_n)\otimes \cF)})\oplus \{0\}\right)]\\
&\ominus\left\{ \Theta_2 f\oplus X_\Theta^{*}(\Delta_2 f\oplus 0):\ f\in F^2(H_n)\otimes \cF\right\}
\end{split}
\end{equation}
and
\begin{equation}\label{H1}
\begin{split}
{\bf H}_1=\{\Theta_2 f\oplus X_\Theta^{*}(\Delta_2f\oplus g)&:\ f\in F^2(H_n)\otimes \cF,\  g\in 
\overline{\Delta_1 (F^2(H_n)\otimes \cE)}\}\\
&\ominus
\left\{ \Theta \varphi\oplus \Delta_\Theta \varphi:\ \varphi\in F^2(H_n)\otimes \cE\right\},
\end{split}
\end{equation}
respectively, where $\Delta_j:=(I-\Theta_j^* \Theta_J)^{1/2}$, $j=1,2$.. Moreover, under the same identification, 
the row contraction $\widetilde{\bf T}$ is unitarily equivalent to the row contraction 
${\bf T}:=[{\bf T}_1,\ldots, {\bf T}_n]$   defined on the Hilbert space
\begin{equation*}
{\bf H}:=[ (F^2(H_n)\otimes {\E_*})\oplus\overline{\Delta_\Theta (F^2(H_n)\otimes \E)}]\ominus\{\Theta g\oplus\Delta_\Theta g :\quad
g\in F^2(H_n)\otimes \cE\},
\end{equation*}
by
$$
{\bf T}_i^*(f\oplus \Delta_\Theta g):=(S_i^*\otimes I_{\cE_*}) f \oplus C_i^*(\Delta_\Theta g),\qquad  i=1,\ldots,n,
$$
where each operator $~C_i~$ is defined by $$~C_i(\Delta_\Theta g):=\Delta_\Theta (S_i\otimes I_\cE) g,\quad 
\ g\in F^2(H_n)\otimes \cE,
$$
 and   $S_1,\ldots, S_n$ are the left creation operators on $F^2(H_n)$.

Since the factorization $\Theta=\Theta_2\Theta_1$  is regular,
 $X_\Theta$ is a unitary operator  which identifies
the subspace
$\overline{\Delta_\Theta(F^2(H_n)\otimes \cE)}$ with
$\overline{\Delta_2(F^2(H_n)\otimes \cF)} \oplus 
\overline{\Delta_1(F^2(H_n)\otimes \cE)}$ and the operator
$C_i$  with $\left(\begin{matrix}
F_i&0\\0&E_i
\end{matrix}\right)$, for  each $i=1,\ldots,n$. Under this identification  the Hilbert spaces ${\bf H}$,   ${\bf H_1}$,
and ${\bf H_2}$ are identified with 
$\HH, \HH_1$, and $\HH_2$, respectively,  and the row contraction ${\bf T}$ is unitarily  equivalent to the row contraction $\TT$.

{\bf Part III.}  We prove the converse of the theorem. 
Due to the above identification, it is enough to assume that the factorization $\Theta=\Theta_2\Theta_1$ is regular  and the subspaces ${\bf H}_1$ and ${\bf H}_2$ are defined as above by relations \eqref{H1} and \eqref{H2}, respectively.
Since $X_\Theta$ is  a unitary operator and using the definition \eqref{X-reg}, we have
\begin{equation*}
\begin{split}
{\bf G}_2:&=
\left\{\Theta_2f\oplus X_\Theta^*(\Delta_2 f\oplus g):\ f\in F^2(H_n)\otimes \cF, \ g\in \overline{\Delta_1 (F^2(H_n)\otimes \cE)}\right\}\\
&\supset
\left\{ \Theta_2 \Theta_1 \varphi \oplus X_\Theta^*(\Delta_2\Theta_1 \varphi \oplus \Delta_1 \varphi):\ \varphi\in F^2(H_n)\otimes \cE\right\}\\
&=
\left\{ \Theta \varphi+ \Delta_\Theta \varphi:\ 
\varphi \in F^2(H_n)\otimes \cE\right\}
\end{split}
\end{equation*}
Hence, we obtain
$$
{\bf H}_1={\bf G}_2\ominus \left\{ \Theta \varphi+ \Delta_\Theta \varphi:\ 
\varphi \in F^2(H_n)\otimes \cE\right\}.
$$
On the other hand,  we have
 \begin{equation*}
\begin{split}
\left[ (F^2(H_n)\right.&\left.\otimes \cE_*)\oplus 
\overline{\Delta_\Theta (F^2(H_n)\otimes \cE)}\right]\ominus {\bf G}_2\\
&=
\left[ (F^2(H_n)\otimes \cE_*)\oplus X_\Theta^*\left(\overline{\Delta_2 (F^2(H_n)\otimes \cF)}\oplus \overline{\Delta_1 (F^2(H_n)\otimes \cE)}\right)\right]\ominus {\bf G}_2\\
&=
\left[ (F^2(H_n)\otimes \cE_*)\oplus X_\Theta^*\left(\overline{\Delta_2 (F^2(H_n)\otimes \cF)}\oplus \{0\}\right)\right]\\
&\phantom{xxxxxxxxxxxxxx}\ominus
\left\{ \Theta_2f\oplus X_\Theta^*(\Delta_2f\oplus \{0\}):
f\in F^2(H_n)\otimes \cF\right\}.
\end{split}
\end{equation*}
Consequently,
$$
{\bf H}_2=\left[ (F^2(H_n)\otimes \cE_*)\oplus 
\overline{\Delta_\Theta (F^2(H_n)\otimes \cE)}\right]\ominus {\bf G}_2.
$$
Hence, and taking into account the definition of ${\bf H}_1$, we deduce that ${\bf H}={\bf H}_1\oplus {\bf H}_2$.

It remains to prove that the subspace ${\bf H}_2$ is 
invariant under the operators ${\bf T}_1^*,\ldots, {\bf T}_n^*$.
If $f\in F^2(H_n)\otimes \cE_*$ and $g\in \overline{\Delta_2 (F^2(H_n)\otimes \cF)}$, then the vector 
$x:=f\oplus X_\Theta^*(g\oplus 0)$ is in ${\bf H}_2$  if and only if
\begin{equation}
\label{eq-def}
\Theta_2^* f+\Delta_2g=0.
\end{equation}
Indeed, using relation \eqref{H2},  one can prove that the condition
$$
\left< f\oplus X_\Theta^*(g\oplus 0), \Theta_2\varphi \oplus 
X_\Theta^*(\Delta_2\varphi\oplus 0)\right>=0\quad \text{ for any }\quad  \varphi\in F^2(H_n)\otimes \cF
$$
 is quivalent to
 \eqref{eq-def}.
Since
\begin{equation*}
\begin{split}
{\bf T}_i^* x&={\bf T}^*_i(f\oplus X_\Theta^*(g\oplus 0))\\
&=(S_i^*\otimes I_{\cE_*})f\oplus C_i^*X_\Theta^*(g\oplus 0)
\end{split}
\end{equation*}
 for each $i=1,\ldots,n$,  to prove that ${\bf T}_i^* x\in {\bf H}_2$,
 it is enough to show that
$$
\left< 
(S_i^*\otimes I_{\cE_*})f\oplus C_i^*(X_\Theta^*(g\oplus 0)),
\Theta_2\varphi \oplus 
X_\Theta^*(\Delta_2\varphi\oplus 0)\right>=0
$$
for any $\varphi\in F^2(H_n)\otimes \cF$.
Since $\Theta$ is a multi-analytic operator, the latter condition is equivalent to 
\begin{equation}\label{SD}
(S_i^*\otimes I_\cF)\Theta_2^*f+\Delta_2 P_1 X_\Theta C_i^*X_\Theta^*(g\oplus 0)=0,
\end{equation}
where $P_1$ is the orthogonal projection of  the direct sum 
$\overline{\Delta_2 (F^2(H_n)\otimes \cF)}\oplus \overline{\Delta_1 (F^2(H_n)\otimes \cE)}$ onto 
$\overline{\Delta_2 (F^2(H_n)\otimes \cF)}$.
Using Lemma \ref{intert-AB} and the definition of the operators $C_i$, $E_i$, and $F_i$, we deduce that
\begin{equation*}
\begin{split}
\Delta_2 P_1 X_\Theta C_i^*X_\Theta^*(g\oplus 0)&=
\Delta_2 P_1 X_\Theta X_\Theta^*\left(\begin{matrix}
F_i^*&0\\
0& E_i^*\end{matrix}
\right) (g\oplus 0)\\
&=\Delta_2F_i^*g\\
&=(S_i^*\otimes I_\cF)\Delta_2 g.
\end{split}
\end{equation*}
Hence, and using relation
\eqref{eq-def}, we have
\begin{equation*}
(S_i^*\otimes I_\cF)\Theta_2^*f+\Delta_2 P_1 X_\Theta C_i^*X_\Theta^*(g\oplus 0)=(S_i^*\otimes I_\cF)(\Theta_2^*f+\Delta_2g)=0,
\end{equation*}
which proves relation \eqref{SD}.
This shows that ${\bf T}_i^*{\bf H}_2\subset {\bf H}_2$ for any $i=1,\ldots,n$. Consequently,  the subspace ${\bf H}_1={\bf H}\ominus {\bf H}_2$
is invariant under the  operators ${\bf T}_1,\ldots,{\bf T}_n$.
This completes the proof of the theorem.
 \end{proof}

  Now we  can reformulate Theorem \ref{inv-factor1} 
in terms of the functional model of a c.n.c. row contraction provided by Theorem \ref{funct-model2}.
This version will be useful later on.

\begin{theorem}\label{inv-factor}
Let $~\Theta:F^2(H_n)\otimes \cE\to F^2(H_n)\otimes \cE_*$ be a purely contractive multi-analytic operator
such that
$$
\overline{\Delta_\Theta (F^2(H_n)\otimes \cE)}=\overline{\Delta_\Theta[ (F^2(H_n)\otimes \cE)\ominus\cE]}.
$$
and let 
 ${\bf T}:=[{\bf T}_1,\ldots, {\bf T}_n]$  be defined on the Hilbert space
\begin{equation*}
{\bf H}:=[ (F^2(H_n)\otimes {\E_*})\oplus\overline{\Delta_\Theta (F^2(H_n)\otimes \E)}]\ominus\{\Theta g\oplus\Delta_\Theta g :\quad
g\in F^2(H_n)\otimes \cE\},
\end{equation*}
by
$$
{\bf T}_i^*(f\oplus \Delta_\Theta g):=(S_i^*\otimes I_{\cE_*}) f \oplus C_i^*(\Delta_\Theta g),\qquad  i=1,\ldots,n,
$$
where each operator $~C_i~$ is defined by $$~C_i(\Delta_\Theta g):=\Delta_\Theta (S_i\otimes I_\cE) g,\quad 
\ g\in F^2(H_n)\otimes \cE,
$$
 and   $S_1,\ldots, S_n$ are the left creation operators on $F^2(H_n)$.

If ~${\bf H}_1\subseteq {\bf H}$~ is an invariant subspace under each operator ${\bf T}_i$,  $i=1,\ldots, n$, then there is a regular  factorization
$$
\Theta=\Theta_2\Theta_1,
$$
where  $\Theta_1:F^2(H_n)\otimes\cE\to F^2(H_n)\otimes\cF$ and 
$\Theta_2:F^2(H_n)\otimes\cF\to F^2(H_n)\otimes\cE_*$
are contractive  multi-analytic operators such that, if $X_\Theta$ is the operator defined by \eqref{X-reg}, then the subspaces ${\bf H}_1$ and ${\bf H}_2:={\bf H}\ominus {\bf H}_1$ have the representations:

  \begin{equation*}
\begin{split}
{\bf H}_1=\{\Theta_2 f\oplus X_\Theta^{*}(\Delta_2f\oplus g)&:\ f\in F^2(H_n)\otimes \cF, \ g\in 
\overline{\Delta_1 (F^2(H_n)\otimes \cE)}\}\\
&\ominus
\left\{ \Theta \varphi\oplus \Delta_\Theta \varphi:\ \varphi\in F^2(H_n)\otimes \cE\right\}
\end{split}
\end{equation*}
and
\begin{equation*}
\begin{split}
{\bf H}_2=[(F^2(H_n)\otimes \cE_*)&\oplus X_\Theta^{*}
\left(
\overline{\Delta_2 (F^2(H_n)\otimes \cF)})\oplus \{0\}\right)]\\
&\ominus\left\{ \Theta_2 f\oplus X_\Theta^{*}(\Delta_2 f\oplus 0):\ f\in F^2(H_n)\otimes \cF\right\}.
\end{split}
\end{equation*}

Conversely,   every regular factorization  $\Theta=\Theta_2\Theta_1$ generates  via the above formulas  the subspaces ${\bf H}_1$  and ${\bf H}_2$  with the following properties:
\begin{enumerate}
\item[(i)]
${\bf H}_1$ is an invariant subspace under each operator ${\bf T}_i$,  $i=1,\ldots, n$;
\item[(ii)]
${\bf H}_2={\bf H}\ominus{\bf H}_1$.
\end{enumerate}
\end{theorem}

\bigskip

 In what follows we need the following factorization result
  for contractive multi-analytic operators \cite{Po-entropy}.
 \begin{lemma}\label{pure}
 Let $\Theta\in  R_n^\infty\bar \otimes B(\cE, \cE')$ be a contractive 
 multi-analytic operator. Then $\Theta$ admits  a unique decomposition 
 $\Theta =\Phi\oplus \Lambda$ with the following properties:
 \begin{enumerate}
 \item[(i)] $\Psi\in R_n^\infty\bar \otimes B(\cE_0, \cE_0')$ is purely contractive, i.e.,
 $\|P_{\cE_0'} \Psi h\|< \|h\|$~ for any $h\in \cE_0$, ~$h\neq 0$;
 \item[(ii)] $\Lambda=I\otimes U\in  R_n^\infty\bar \otimes B(\cE_u, \cE_u')$,
  where $U\in B(\cE_u, \cE_u')$ is a unitary operator;
 \item[(iii)]
 $\cE=\cE_0\oplus \cE_u$ and $\cE'=\cE_0'\oplus \cE_u'$.
 \end{enumerate}
Moreover, 
 the purely contractive part of an outer or inner multi-analytic operator is also outer or inner, respectively.
\end{lemma}

 The next result is an addition to Theorem \ref{funct-model2}.

\begin{proposition}\label{pure-cnc}
Let $~\Theta:F^2(H_n)\otimes \cE\to F^2(H_n)\otimes \cE_*$ be a  contractive multi-analytic operator
such that
$$
\overline{\Delta_\Theta (F^2(H_n)\otimes \cE)}=\overline{\Delta_\Theta[ (F^2(H_n)\otimes \cE)\ominus\cE]}.
$$
and let 
 ${\bf T}:=[{\bf T}_1,\ldots, {\bf T}_n]$  be the functional model associated with $\Theta$, as in Theorem $\ref{funct-model2}$.
\begin{enumerate}
\item[(i)]
The characteristic function of ${\bf T}:=[{\bf T}_1,\ldots, {\bf T}_n]$ coincides with the purely contractive part
of $\Theta$.
\item[(ii)]
The space ${\bf H}$ defined in Theorem $\ref{funct-model2}$ is different from
$\{0\}$ if and only if there is no unitary operator $U\in B(\cE,\E_*)$ such that $\Theta=I\otimes U$.
\end{enumerate} 
\end{proposition}

\begin{proof}
According to Lemma \ref{pure}, the multi-analytic operator
$\Theta$ admits the decomposition $\Theta =\Phi\oplus \Lambda$   with $\Psi\in R_n^\infty\bar \otimes B(\cE_0, \cE_{*0})$  purely contractive and 
$\Lambda=I\otimes U\in  R_n^\infty\bar \otimes B(\cE_u, \cE_{*u})$,
  where $U\in B(\cE_u, \cE_{*u})$ is a unitary operator,
 $\cE=\cE_0\oplus \cE_u$,  and $\cE_*=\cE_{*0}\oplus \cE_{*u}$.
Notice that 
\begin{equation*}
\begin{split}
F^2(H_n)\otimes \cE_*&=(F^2(H_n)\otimes \cE_{*u})\oplus (F^2(H_n)\otimes \cE_{*0}) \text{ and}\\
F^2(H_n)\otimes \cE&=(F^2(H_n)\otimes \cE_{u})\oplus (F^2(H_n)\otimes \cE_{0}).
\end{split}
\end{equation*}
On the other hand, we have
\begin{equation*}
\left\{ \Theta g\oplus \Delta_\Theta g:\ g\in F^2(H_n)\otimes \cE\right\}
=
(F^2(H_n)\otimes \cE_{*u})\oplus \left\{\Phi \varphi \oplus \Delta_\Phi \varphi: \ \varphi\in F^2(H_n)\otimes \cE_{0}\right\}.
\end{equation*}
Now, using the definition of the Hilbert space ${\bf H}$, one can identify ${\bf H}$ with 
$$
{\bf H}_0:=\left[ (F^2(H_n)\otimes \cE_{*0})\oplus \overline{\Delta_\Phi(F^2(H_n)\otimes \cE_{0})}\right]\ominus
\left\{\Phi \varphi \oplus \Delta_\Phi \varphi: \ \varphi\in F^2(H_n)\otimes \cE_{0}\right\}.
$$
Due to this identification, the row contraction ${\bf T}:=[{\bf T}_1,\ldots, {\bf T}_n]$ is unitarily equivalent to
${\bf T}^0:=[{\bf T}_1^0,\ldots, {\bf T}_n^0]$, which is defined on ${\bf H}_0$   in the same manner  as ${\bf T}$ is 
 defined on ${\bf H}$.
Since $\Delta_\Theta=\Delta_\Phi\oplus 0$, it is easy to see that
$$
\overline{\Delta_\Phi (F^2(H_n)\otimes \cE)}=\overline{\Delta_\Phi[ (F^2(H_n)\otimes \cE)\ominus\cE]}.
$$
According to the second part of Theorem \ref{funct-model2}
the characteristic function  of ${\bf T}^0$ coincides with 
the multi-analytic operator $\Phi$ which coincides with 
the characteristic function  of ${\bf T}$.

 We prove now part (ii).  If $\Theta=I\otimes U$ for some unitary operator $U\in B(\cE,\E_*)$, then $\Delta_\Theta=0$ and 
$$
{\bf H}=[F^2(H_n)\otimes \cE_*]\ominus \Theta(F^2(H_n)\otimes \cE)=\{0\}.
$$

If $\Theta$ is not a unitary multi-analytic operator, then, according to Lemma \ref{pure}, it has a non-trivial purely contractive part.
By part (i), Theorem \ref{funct-model1}, and Theorem \ref{funct-model2}, we deduce that
$$
\dim \cD_*=\dim\cE_{*0},\quad \dim \cD=\dim\cE_{0},
$$
where $\cE$ and $\cE_{*0}$ are not both equal to $\{0\}$.
Since $\cD_*\subset \cH$ and $\cD\subset \cH^{(n)}$, we deduce that $\cH\neq \{0\}$.
This completes the proof.
\end{proof}

The following result  is   an important addition to Theorem \ref{inv-factor} (and hence  also  to Theorem \ref{inv-factor1}).

\begin{theorem}
\label{inv-factor2}
 Under the conditions of Theorem $\ref{inv-factor}$,
let  ${\bf H}={\bf H}_1\oplus {\bf H}_2$ be the decomposition corresponding to the regular factorization 
$\Theta=\Theta_2\Theta_1$, and let 
$$
{\bf T}_i=\left( \begin{matrix}
{\bf A}_i& *\\
0& {\bf B}_i
\end{matrix}\right), \quad i=1,\ldots, n,
$$
be the corresponding triangulation of ${\bf T}:=[{\bf T}_1,\ldots, {\bf T}_n]$.
Then the characteristic functions of  the row contractions ${\bf A}:=[{\bf A}_1,\ldots, {\bf A}_n]$ and ${\bf B}:=[{\bf B}_1,\ldots, {\bf B}_n]$
coincide with the purely contractive parts of the multi-analytic operators $\Theta_1$ and $\Theta_2$, respectively.

Moreover,
the invariant subspace ${\bf H}_1$ under  the operators ${\bf T}_1,\ldots, {\bf T}_n$ is non-trivial if and only if  the regular factorization 
$\Theta=\Theta_2\Theta_1$  is non-trivial, i.e., each factor is not a unitary constant.
\end{theorem}

\begin{proof}
Define the operator $U$ from the Hilbert space 
$$
 (F^2(H_n)\otimes \cE_*)\oplus X_\Theta^*\left(\overline{\Delta_2 (F^2(H_n)\otimes \cF)}\oplus \{0\}\right)
$$
to
$$
(F^2(H_n)\otimes \cE_*)\oplus \overline{\Delta_2 (F^2(H_n)\otimes \cF)} 
$$
by setting
$$
U(f\oplus X^*(g\oplus 0)):=f\oplus g,
$$
for any $f\in F^2(H_n)\otimes \cE_*$ and $g\in \overline{\Delta_2 (F^2(H_n)\otimes \cF)}$.
Since $X_\Theta$ is unitary, so is $U$. Using the definition of ${\bf H}_2$
(see relation \eqref{H2}), we deduce that
$U{\bf H}_2=\widehat{\cH}_2$, where
\begin{equation}\label{H22}
\begin{split}
\widehat{\cH}_2:=\left[ (F^2(H_n)\otimes \cE_*)\right.&\left.\oplus \overline{\Delta_2 (F^2(H_n)\otimes \cF)}\right]\\
&\ominus\left\{\Theta_2\varphi \oplus \Delta_2 \varphi:\ \varphi\in
F^2(H_n)\otimes \cF\right\}.
\end{split}
\end{equation}
Set $\Gamma_i^*:= U{\bf B}^*_i U^*$, $i=1,\ldots, n$, and 
denote  by $P_1$ the orthogonal projection of the direct sum
$
 \overline{\Delta_2 (F^2(H_n)\otimes \cF)}\oplus \overline{\Delta_1 (F^2(H_n)\otimes \cE)}
$
onto 
$\overline{\Delta_2 (F^2(H_n)\otimes \cF)}$.
Using Lemma \ref{intert-AB}, we deduce that
\begin{equation*}
\begin{split}
P_1X_\Theta C_i^*X_\Theta^*(g\oplus 0)&=
P_1\left(\begin{matrix}F_i^*&0\\0&E_i^*
\end{matrix}\right)
\left(\begin{matrix} g\\0\end{matrix}\right)\\
&=F_i^*g
\end{split}
\end{equation*}
for any $g\in \overline{\Delta_2 (F^2(H_n)\otimes \cF)}$
and $i=1,\ldots, n$.
Hence and using the definitions for the row contraction
$[{\bf T}_1,\ldots, {\bf T}_n]$ and the unitary operator $U$, we have
\begin{equation*}
\begin{split}
\Gamma_i^*(f\oplus g)&=
U{\bf T}_i^*(f\oplus X_\Theta^*(g\oplus 0))\\
&=U\left[(S_i^*\otimes I_{\cE_*})f\oplus C_i^*X_\Theta^*(g\oplus 0)\right]\\
&=(S_i^*\otimes I_{\cE_*})f\oplus P_1 X_\Theta C_i^*X_\Theta^*(g\oplus 0)\\
&=(S_i^*\otimes I_{\cE_*})f\oplus F_i^*g
\end{split}
\end{equation*}
for any $f\in F^2(H_n)\otimes \cE_*$ and $g\in \overline{\Delta_2 (F^2(H_n)\otimes \cF)}$ such that
$f\oplus g\in \cH_2$, and $i=1,\ldots, n$.

Since 
$$
\overline{\Delta_\Theta(F^2(H_n)\otimes \cE)}=
\overline{\Delta_\Theta(F^2(H_n)\otimes \cE)\ominus\cE},
$$
one can use  again Lemma \ref{intert-AB} to deduce that
$$
\overline{\Delta_2(F^2(H_n)\otimes \cF)}=
\overline{\Delta_2(F^2(H_n)\otimes \cF)\ominus\cF}.
$$
Now, due to  Proposition \ref{pure-cnc}, we infer that the characteristic function of 
the row contraction $[\Gamma_1,\ldots, \Gamma_n]$, $\Gamma_i\in B(\widehat{\cH}_2)$, (and hence 
also $[{\bf B}_1,\ldots, {\bf B}_n]$) coincides with the purely contractive part of the multi-analytic operator $\Theta_2$.

Taking into account the definition of the subspace ${\bf H}_1$
(see relation \eqref{H1}) and the fact that $\Theta=\Theta_2 \Theta_1$, one can see that, for each 
$f\in F^2(H_n)\otimes \cF$ and $g\in \overline{\Delta_1 (F^2(H_n)\otimes \cE)}$, the vector 
$\Theta_2 f\oplus X_\Theta^*(\Delta_2 f\oplus g)$ is in ${\bf H}_1$ if and only if 
$$
\left<\Theta_2 f\oplus X_\Theta^*(\Delta_2 f\oplus g), \Theta_2 \Theta_1\varphi \oplus X_\Theta^*(\Delta_2 \Theta_1 \varphi \oplus \Delta_1\varphi)\right>=0
$$
for any $\varphi\in F^2(H_n)\otimes \cE$.
The latter  equation is equivalent to
$$
\Theta_1^*\Theta_2^* \Theta_2f+\Theta_1^* \Delta_2^2f+\Delta_1g=0.
$$
Since $\Delta_2^2=I-\Theta_2^*\Theta_2$, the above equation is
equivalent to 
\begin{equation}
\label{Th-De}
\Theta_1^*f +\Delta_1g=0.
\end{equation}
If $x:=\Theta_2 f\oplus X_\Theta^*(\Delta_2 f\oplus g)\in {\bf H}_1$, then we have
$$
{\bf T}_i^*x=(S_i^*\otimes I_{\cE_*})\Theta_2f\oplus C_i^*X_\Theta^*(\Delta_2 f\oplus g)
$$
for each $i=1,\ldots, n$.
Since $\Theta_2$ is a multi-analytic operator
and 
$$
f=\sum_{j=1}^n (S_j  S_j^*\otimes I_\cF)f+ f(0),
$$
where $f(0):=P_{1\otimes \cF} f$, we deduce that
\begin{equation*}
\begin{split}
{\bf T}_i^*x&=[\Theta_2(S_i^*\otimes I_{\cF})f+(S_i^*\otimes I_{\cE_*})\Theta_2f(0)]\oplus C_i^*X_\Theta^*(\Delta_2 f\oplus g)\\
&=u+v,
\end{split}
\end{equation*}
where
$$
u:=\Theta_2(S_i^*\otimes I_{\cF})f\oplus \left[X_\Theta^*(\Delta_2 (S_i^*\otimes I_\cF)f\oplus E_i^*g)\right]
$$
and 
$$
v:=(S_i^*\otimes I_{\cE_*})\Theta_2 f(0)\oplus\left[
C_i^*X_\Theta^*(\Delta_2 f\oplus g)-X_\Theta^*(\Delta_2 (S_i^*\otimes I_\cF)f\oplus E_i^*g)\right].
$$
Now notice that
$u\in {\bf H}_1$. Indeed, using the above characterization of the 
elements of ${\bf H}_1$, it is enough to show that
\begin{equation}
\label{Th*}
\Theta_1^* (S_i^*\otimes  I_{\cF}) f+\Delta_1 E_i^* g=0,\quad i=1,\ldots,n.
\end{equation}
Using relation \eqref{Th-De} and the definition
of $E_i$, we have
\begin{equation*}
\begin{split}
\Theta_1^* (S_i^*\otimes  I_{\cF}) f+\Delta_1 E_i^* g&=
(S_i^*\otimes I_\cE)(\Theta_1^* f+ \Delta_1 g)\\
&=0,
\end{split}
\end{equation*}
which proves \eqref{Th*} and therefore $u\in {\bf H}_1$.

Now we prove that $v\in {\bf H}_2$. First, notice that due to Lemma \ref{intert-AB}, we have 
$$C_i^* X_\Theta^*(0\oplus g)=X_\Theta^*(0\oplus E_i^*g),\quad
g\in \overline{\Delta_1(F^2(H_n)\otimes \cE)},
$$
and therefore
\begin{equation}\label{v}
v=(S_i^*\otimes I_{\cE_*})\Theta_2 f(0)\oplus\left[
C_i^*X_\Theta^*(\Delta_2 f\oplus 0)-X_\Theta^*(\Delta_2 (S_i^*\otimes I_\cF)f\oplus 0)\right].
\end{equation}
Using again Lemma \ref{intert-AB}, and the definition of $F_i$, we infer that
\begin{equation*}
\begin{split}
C_i^*X_\Theta^*(\Delta_2 f\oplus 0)&=
C_i^*X_\Theta^*\left(\Delta_2\left(\sum_{j=1}^n S_jS_j^*\otimes I_\cF\right) f(0)\oplus 0\right)
+C_i^*X_\Theta^*(\Delta_2 f(0)\oplus 0)\\
&=
X_\Theta^*\left( F_i^* \Delta_2\left(\sum_{j=1}^n S_jS_j^*\otimes I_\cF\right) f\oplus 0\right)+C_i^*X_\Theta^*(\Delta_2 f(0)\oplus 0)\\
&=
X_\Theta^*(\Delta_2(S_i^*\otimes I_\cF)f\oplus 0)+C_i^*X_\Theta^*(\Delta_2 f(0)\oplus 0)\\
&=
X_\Theta^*(\Delta_2(S_i^*\otimes I_\cF)f\oplus 0)+
X_\Theta^*(F_i^*\Delta_2 f(0)\oplus 0).
\end{split}
\end{equation*}
Consequently, relation \eqref{v} implies
$$
v=(S_i^*\otimes I_{\cE_*})\Theta_2 f(0)\oplus X_\Theta^*(F_i^*\Delta_2 f(0)\oplus 0).
$$
Due to the definition of the subspace ${\bf H}_2$, to prove that $v\in {\bf H}_2$, it is enough to show that
$$
\Theta_2^*(S_i^*\otimes I_{\cE_*}) \Theta_2 f(0)+\Delta_2 F_i^* \Delta_2 f(0)=0
$$
for each $i=1,\ldots, n$.
Since 
$$
\Delta_2  F_i^*=(S_i^*\otimes I_\cF) \Delta_2,\quad i=1,\ldots, n, 
$$
 and $\Theta_2$ is multi-analytic, we have
\begin{equation*}
\begin{split}
\Theta_2^*(S_i^*\otimes I_{\cE_*}) \Theta_2 f(0)+\Delta_2 F_i^* \Delta_2 f(0)&=
(S_i^*\otimes I_\cF)(\Theta_2^* \Theta_2+\Delta_2^2) f(0)\\
&=(S_i^*\otimes I_\cF)  f(0)=0.
\end{split}
\end{equation*}
Hence, $v\in {\bf H}_2$. Now, using the  fact that ${\bf T}_i^* x=u+v$ and the definitions for $u$ and $v$, we deduce  that the operator
${\bf A}^*_i:=P_{{\bf H}_1}{\bf T}_i^*|_{{\bf H}_1}$ satisfies the 
equation
\begin{equation}
\label{A1}
{\bf A}_i^*\left(\Theta_2f\oplus X_\Theta^*(\Delta_2 f\oplus g)\right)=\Theta_2(S_i^*\otimes I_{\cF})f\oplus \left[X_\Theta^*(\Delta_2 (S_i^*\otimes I_\cF)f\oplus E_i^*g)\right]
\end{equation}
for any  $\Theta_2f\oplus X_\Theta^*(\Delta_2 f\oplus g)\in {\bf H}_1$ and  $i=1,\ldots, n$.

Now, define the operator $\Omega$ from
$$
\left\{ \Theta_2 f\oplus X_\Theta^*(\Delta_2 f\oplus g):\ f\in F^2(H_n)\otimes \cF, g\in 
\overline{\Delta_1(F^2(H_n)\otimes \cE)}\right\}
$$
to the direct sum \ 
$(F^2(H_n)\otimes \cF)\oplus \overline{\Delta_1(F^2(H_n)\otimes \cE)}
$~
by setting
\begin{equation} \label{ome}
\Omega(\Theta_2 f\oplus X_\Theta^*(\Delta_2f\oplus g)):=f\oplus g.
\end{equation}
Since 
\begin{equation*}
\begin{split}
\|\Theta_2f\oplus X_\Theta^*(\Delta_2f\oplus g)\|^2&=
\|\Theta_2f\|^2+\|X_\Theta^*(\Delta_2f\oplus g)\|^2\\
&=
\left<\Theta_2^* \Theta_2f,f\right>+\|\Delta_2f\|^2+\|g\|^2\\
&=\|f\oplus g\|^2,
\end{split}
\end{equation*}
it is clear that $\Omega$ is a unitary operator. Notice also that
\begin{equation*}
\begin{split}
\Omega(\Theta\varphi\oplus \Delta_\Theta \varphi)&=
\Omega(\Theta_2 \Theta_1\varphi\oplus X_\Theta^*(\Delta_2 \Theta_1 \varphi\oplus \Delta_1\varphi))\\
&=\Theta_1 \varphi\oplus \Delta_1\varphi
\end{split}
\end{equation*}
 for any $\varphi\in F^2(H_n)\otimes \cE$.
Consequently, $\Omega {\bf H}_1=\widehat{\cH}_1$, where
\begin{equation}\label{H11}
\begin{split}
\widehat{\cH}_1:=[ (F^2(H_n)\otimes \cF)&\oplus \overline{\Delta_1(F^2(H_n)\otimes \cE)}]\\
&\ominus\left\{ \Theta_1 \varphi\oplus \Delta_1\varphi:\ \varphi\in F^2(H_n)\otimes \cE\right\}.
\end{split}
\end{equation}
Setting $\Lambda_i:=\Omega {\bf A_i} \Omega^*$, \  relation \eqref{A1} implies
$$
\Lambda_i^*(f\oplus g)=(S_i^*\otimes I_\cF)f \oplus E_i^*g, \quad  f\oplus g\in \cH_1,
$$
 for any $i=1,\ldots,n$.
Once again, 
Lemma \ref{intert-AB}  implies
$$
\overline{\Delta_1(F^2(H_n)\otimes \cE)}=
\overline{\Delta_1(F^2(H_n)\otimes \cE)\ominus\cE}.
$$
Now, using Proposition \ref{pure-cnc}, we infer that the characteristic function of 
the row contraction $[\Lambda_1,\ldots, \Lambda_n]$, $\Lambda_i\in B(\widehat{\cH}_1)$,  (and hence also
$[{\bf A}_1,\ldots, {\bf A}_n]$) coincides with the purely contractive part of the multi-analytic operator $\Theta_1$.
Due to the relations \eqref{H22}, \eqref{H11}, and Proposition \ref{pure-cnc}, the subspaces $\widehat{\cH}_1$ and $\widehat{\cH}_2$ (and hence also ${\bf H}_1$ and ${\bf H}_2$) are different from $\{0\}$
if and only if both multi-analytic operators $\Theta_1$ and $\Theta_2$ are not unitary constant, i.e., the factorization
$\Theta=\Theta_1 \Theta_2$ is non-trivial.
This completes the proof.
 \end{proof}


Now, combining  Theorem \ref{inv-factor1} and Theorem \ref{inv-factor2},  we can deduce the following result.

\begin{theorem}
\label{non-trivial sub}
Let $T:=[T_1,\ldots, T_n]$ be a completely non-coisometric row contraction
on a separable Hilbert space $\cH$. Then, there is a non-trivial invariant subspace under each operator $T_1,\ldots, T_n$ 
if and only if the characteristic function $\Theta_T$  has a non-trivial regular factorization.
\end{theorem}

Concerning the uniqueness in   Theorem \ref{inv-factor} (and also Theorem \ref{inv-factor1}), we can prove the following result, which shows the extent
to which a joint invariant subspace determines the corresponding regular factorization
of the characteristic function.

\begin{theorem}
\label{inv-factor3}
 Under the conditions of Theorem $\ref{inv-factor}$,
let
$$
\Theta=\Theta_2\Theta_1\quad \text{ and }\quad \Theta=\Theta_2'\Theta_1'
$$ be  two regular factorizations of the   purely contractive multi-analytic operator $\Theta$, and let $\cE,\cF,\cE_*$, and $\cE,\cF',\cE_*$ be the corresponding  Hilbert spaces.
Let ${\bf H}_1\subset {\bf H}$ and ${\bf H}_1'\subset {\bf H}$
be the invariant subspaces under each operator ${\bf T}_i$,  $i=1,\ldots, n$, corresponding to the above factorizations.
If 
${\bf H}_1\subset {\bf H}_1'$, then there is a multi-analytic operator
$\Psi:F^2(H_n)\otimes \cF\to  F^2(H_n)\otimes \cF'$ such that
$$
\Theta_1'=\Psi \Theta_1.
$$

Moreover, if
${\bf H}_1= {\bf H}_1'$, then 
$$
\Theta_1'=(I\otimes \Psi_0) \Theta_1
$$
for some unitary operator $\Psi_0\in B(\cF,\cF')$ and, consequently,  the multi-analytic operators $\Theta_1$ and $\Theta_1'$ coincide.
\end{theorem}

\begin{proof}
We associate with the factorization $\Theta=\Theta_2\Theta_1$ the subspace
$$
\cM:=\left\{\Theta_2f\oplus X_\Theta^*(\Delta_2 f\oplus g):\ f\in F^2(H_n)\otimes \cF, g\in \overline{\Delta_1(F^2(H_n)\otimes \cE)}\right\}.
$$
Similarly, we define the subspace $\cM'$ associated with the factorization $\Theta=\Theta_2'\Theta_1'$. Since ${\bf H}_1\subseteq {\bf H}_1'$, relation \eqref{H1}  and its analogue for ${\bf H}_1'$ imply $\cM\subseteq \cM'$.
Consequently, for each $f\in F^2(H_n)\otimes \cF$, there exist $f'\in F^2(H_n)\otimes \cF'$ and  $g'\in \overline{\Delta_1'(F^2(H_n)\otimes \cE)}
$ such that
\begin{equation}
\label{Th2}
\Theta_2f\oplus X_\Theta^*(\Delta_2 f\oplus 0)=\Theta_2'f'\oplus X_\Theta'^*(\Delta_2 f'\oplus g').
\end{equation}
Hence and using the definition of the unitary operators $X_\Theta$ and $X_\Theta'$, we have 
\begin{equation*}
\begin{split}
\|f\|^2&=\|\Theta_2f\oplus X_\Theta^*(\Delta_2 f\oplus g)\|^2\\
&=\|\Theta_2'f'\oplus X_\Theta'^*(\Delta_2 f'\oplus g')\|^2\\
&=\|f'\|^2+\|g'\|^2.
\end{split}
\end{equation*}
Therefore, it makes sense to define the contractions $Q:F^2(H_n)\otimes \cF\to F^2(H_n)\otimes \cF'$ and 
$R:F^2(H_n)\otimes \cF\to\overline{\Delta_1'(F^2(H_n)\otimes \cE)}$ by setting
$Qf:=f'$ and $Rf:=g'$, respectively.
Now, we show that $Q$ is a multi-analytic operator, i.e.,
$$Q(S_i\otimes I_\cF)=(S_i\otimes I_{\cF'})Q,\quad i=1,\ldots, n.
$$
Let $f_1,\ldots, f_n$ be arbitrary elements in $F^2(H_n)\otimes \cE$.
Taking into account the definitions for $C_i$ and $X_\Theta$, and the fact that
$$
(S_j^*\otimes I_\cF)\Delta_2^2(S_i\otimes I_{\cF})=\delta_{ij} \Delta_2^2,\quad i,j=1,\ldots,n,
$$
we deduce that
\begin{equation*}
\begin{split}
 \left< C_iX_\Theta^*(\Delta_2 f\oplus 0),\Delta_\Theta\left(\sum_{j=1}^n(S_j\otimes I_\cE)f_j\right)\right>
&=
\left<
\Delta_2 f\oplus 0), X_\Theta\Delta_\Theta f_i
\right>\\
&=
\left<
\Delta_2 f\oplus 0),\Delta_2\Theta_1 f_i\oplus \Delta_1f_i
\right>\\
&=\left<\Delta_2^2 f, \Theta_1f_i\right>
\end{split}
\end{equation*}
and
\begin{equation*}
\begin{split}
\Biggl< X_\Theta^*(\Delta_2(S_i\Biggr.&\Biggl.\otimes I_\cF)f\oplus 0), \Delta_\Theta
\left(\sum_{j=1}^n(S_j\otimes I_\cE)f_j\right)\Biggr>\\
&=
\left<\Delta_2(S_i\otimes I_\cF)f\oplus 0,
\Delta_2 \Theta_1 \left(\sum_{j=1}^n(S_j\otimes I_\cE)f_j\right)\oplus \Delta_1 \left(\sum_{j=1}^n(S_j\otimes I_\cE)f_j\right)\right>\\
&=
\left<\Delta_2(S_i\otimes I_\cF)f,
\Delta_2 \Theta_1 \left(\sum_{j=1}^n(S_j\otimes I_\cE)f_j\right)\right>\\
&=\sum_{j=1}^n \left<(S_j^*\otimes I_\cF)\Delta_2^2(S_i\otimes I_{\cF})f, \Theta_1 f_j\right>\\
&=
\left<\Delta_2^2 f, \Theta_1f_i\right>.
\end{split}
\end{equation*}
Hence, and taking into account that
$$
\overline{\Delta_\Theta(F^2(H_n)\otimes \cE)}=
\overline{\Delta_\Theta [(F^2(H_n)\otimes \cE)\ominus \cE]},
$$
we deduce that
\begin{equation}
\label{CiX}
C_iX_\Theta^*(\Delta_2 f\oplus 0)=
X_\Theta^*(\Delta_2(S_i\otimes I_\cF)f\oplus 0)\quad \text{ for any } \quad f\in F^2(H_n)\otimes \cF.
\end{equation}
Similar calculations show that 
\begin{equation}\label{CiX2}
C_iX_\Theta^*(0\oplus \Delta_1\varphi)=
X_\Theta^*(0\oplus \Delta_1(S_i\otimes I_\cE)\varphi)
\end{equation}
for any $\varphi\in F^2(H_n)\otimes \cE$ and $i=1,\ldots,n$.
Moreover,  similar relations to \eqref{CiX}  and \eqref{CiX2} hold with $X_\Theta'$, $\Delta_1'$, and $\Delta_2'$ instead of $X_\Theta$, $\Delta_1$, and $\Delta_2$, respectively.
 Since
\begin{equation}
\label{CX'}
C_iX_\Theta'^*(0\oplus \Delta_1'\varphi)=
X_\Theta'^*(0\oplus \Delta_1'(S_i\otimes I_\cE)\varphi)
\end{equation}
for any $\varphi\in F^2(H_n)\otimes \cE$ and $i=1,\ldots,n$,
by taking appropriate limits, we deduce that
$$
C_iX_\Theta'^*(\{0\}\oplus\overline{\Delta_1'(F^2(H_n)\otimes \cE)})\subseteq X_\Theta'^*(\{0\}\oplus  \overline{\Delta_1'(F^2(H_n)\otimes \cE)}).
$$
Consequently,
for each $g'\in \overline{\Delta_1'(F^2(H_n)\otimes \cE)}$ there exists 
$g''\in \overline{\Delta_1'(F^2(H_n)\otimes \cE)}$ 
such that
\begin{equation}
\label{CiX'}
C_iX_\Theta'^*(0\oplus g')=X_\Theta'^*(0\oplus g'').
\end{equation}
Now, notice that using relations 
\eqref{CiX}, \eqref{Th2}, \eqref{CX'},  and \eqref{CiX'}, we obtain
\begin{equation*}
\begin{split}
\Theta_2(S_i\otimes I_\cF)f\oplus X_\Theta^*(\Delta_2(S_i\otimes I_\cF)f\oplus 0)&=
(S_i\otimes I_{\cE_*}\oplus C_i) (\Theta_2 f\oplus X_\Theta^*(\Delta_2f\oplus 0))\\
&=
(S_i\otimes I_{\cE_*}\oplus C_i) (\Theta_2' f'\oplus X_\Theta'^*(\Delta_2'f'\oplus g'))\\
&=
\Theta_2'(S_i\otimes I_{\cF'})f'\oplus 
X_\Theta'^*(\Delta_2'(S_i\otimes I_{\cF'})f'\oplus g'')
\end{split}
\end{equation*}
for any $f\in F^2(H_n)\otimes \cF$.
Hence and using the definition of $Q$, we deduce that
$$
Q(S_i\otimes I_\cF)f= (S_i\otimes I_{\cF'})f'=(S_i\otimes I_{\cF'})Qf,\quad f\in F^2(H_n)\otimes \cF,
$$
which proves that $Q$ is a multi-analytic operator.

Since $\cM\subset \cM'$, we have
\begin{equation}
\label{inte}
\bigcap_{k=0}^\infty \bigoplus_{|\alpha|=k}[(S_\alpha\otimes I_{\cE_*})\oplus C_\alpha]\cM\subseteq 
\bigcap_{k=0}^\infty \bigoplus_{|\alpha|=k}[(S_\alpha\otimes I_{\cE_*})\oplus C_\alpha]\cM'.
\end{equation}
Using Lemma \ref{intert-AB}, the definition  \eqref{ome} of the unitary operator $\Omega$, and relations \eqref{CiX}, \eqref{CiX2}, one can prove that
$$
[(S_i\otimes I_{\cE_*})\oplus C_i]\Omega^*=\Omega^*[(S_i\otimes I_\cF)\oplus E_i].
$$
Indeed, we have
\begin{equation*}
\begin{split}
[(S_i\otimes I_{\cE_*})\oplus C_i]\Omega^*(f\oplus \Delta_1 \varphi)
&=
\Theta_2(S_i\otimes I_\cF)f\oplus C_iX_\Theta^*(\Delta_2 f\oplus \Delta_1 \varphi)\\
&=
\Theta_2(S_i\otimes I_\cF)f\oplus 
X_\Theta^*(\Delta_2 (S_i\otimes I_\cF)f\oplus \Delta_1 (S_i\otimes I_\cE)\varphi)\\
&=
\Omega^*[(S_i\otimes I_\cF)f\oplus \Delta_1 (S_i\otimes I_\cE)\varphi)]\\
&=
\Omega^*[(S_i\otimes I_\cF)\oplus E_i](f\oplus \Delta_1 \varphi)
\end{split}
\end{equation*}
for any $f\in F^2(H_n)\otimes \cF$ and $\varphi\in F^2(H_n)\otimes \cE$.

Now,  due to  the fact that $[S_1\otimes I_\cF,\ldots, S_n\otimes I_\cF]$ is a multi-shift and $[E_1,\ldots, E_n]$
is a Cuntz row isometry,  the noncommutative Wold decomposition implies
\begin{equation*}
\begin{split}
\bigcap_{k=0}^\infty &\bigoplus_{|\alpha|=k}[(S_\alpha\otimes I_{\cE_*})\oplus C_\alpha]\cM\\
&=\Omega^*\left\{
\bigcap_{k=0}^\infty \left[\bigoplus_{|\alpha|=k}(S_\alpha\otimes I_\cF)(F^2(H_n)\otimes \cF)\right]\oplus 
\bigcap_{k=0}^\infty \left[\bigoplus 
 \limits_{|\alpha|=k} E_\alpha \overline{\Delta_1(F^2(H_n)\otimes \cE)}\right]
\right\}\\
&=\Omega^*\left( \{0\}\oplus \overline{\Delta_1(F^2(H_n)\otimes \cE)}\right)\\
&=
\left\{ 0\oplus X_\Theta^*(0\oplus g):\ g\in \overline{\Delta_1(F^2(H_n)\otimes \cE)}\right\}.
\end{split}
\end{equation*}
A similar relation can be obtain for the set on the  right side of  the inclusion \eqref{inte}. Hence and  using relation \eqref{inte}, we obtain
$$
\left\{ 0\oplus X_\Theta^*(0\oplus g):\ g\in \overline{\Delta_1(F^2(H_n)\otimes \cE)}\right\}
\subseteq 
\left\{ 0\oplus X_\Theta'^*(0\oplus g'):\ g'\in \overline{\Delta_1'(F^2(H_n)\otimes \cE)}\right\}.
$$
Consequently, for each $g\in \overline{\Delta_1(F^2(H_n)\otimes \cE)}$ there exists
$g'\in \overline{\Delta_1'(F^2(H_n)\otimes \cE)}$ such that
\begin{equation}
\label{XX'}
X_\Theta^*(0\oplus g)=X_\Theta'^*(0\oplus g').
\end{equation}
Since $X_\Theta$ and $X_\Theta'$ are unitary operators, we can define the isometry
$$V:\overline{\Delta_1(F^2(H_n)\otimes \cE)}\to
\overline{\Delta_1'(F^2(H_n)\otimes \cE)}
$$
by setting $Vg:=g'$.
For each $\varphi\in F^2(H_n)\otimes \cE$, we have
\begin{equation}
\label{TDT}
\Theta \varphi\oplus \Delta_\Theta \varphi= \Theta_2'\Theta_1' \varphi\oplus
X_\Theta'^*(\Delta_2' \Theta_1' \varphi\oplus \Delta_1' \varphi).
\end{equation}
On the other hand, using the operators $Q,R,V$ and relation 
\eqref{Th2}, we deduce that

\begin{equation*}
\begin{split}
\Theta \varphi\oplus \Delta_\Theta \varphi
&= 
\Theta_2\Theta_1 \varphi\oplus
X_\Theta^*(\Delta_2 \Theta_1 \varphi\oplus \Delta_1 \varphi)\\
&=
\left[\Theta_2\Theta_1 \varphi\oplus
X_\Theta^*(\Delta_2 \Theta_1 \varphi\oplus 0\right]+
\left[0\oplus X_\Theta^*(0\oplus \Delta_1 \varphi)\right]\\
&=
\left[\Theta_2'Q\Theta_1 \varphi \oplus X_\Theta'^*(\Delta_2'Q\Theta_1 \varphi \oplus R\Theta_1 \varphi)\right]
+
\left[ 0\oplus X_\Theta'^*(0\oplus V\Delta_1 \varphi )\right]\\
&=\Theta_2'Q\Theta_1 \varphi \oplus X_\Theta'^*(\Delta_2'Q\Theta_1 \varphi \oplus y),
\end{split}
\end{equation*}
where $y:= R\Theta_1 \varphi+ V\Delta_1 \varphi$ is in 
$\overline{\Delta_1'(F^2(H_n)\otimes \cE)}$.
Using the latter relation and \eqref{TDT}, we obtain
$$
\Theta_2'\Theta_1'\varphi= \Theta_2'Q \Theta_1 \varphi \quad \text{ and } \quad 
\Delta_2'\Theta_1'\varphi=
\Delta_2' Q\Theta_1 \varphi.
$$
Since the mapping $\Theta_2'f'\oplus \Delta_2' f'\mapsto
f'$  is isometric, we deduce that
\begin{equation}
\label{first}
\Theta_1'\varphi =Q\Theta_1 \varphi,\quad \varphi\in F^2(H_n)\otimes \cE,
\end{equation}
which proves the first part of the theorem.

Now assume that ${\bf H}_1={\bf H}_1'$. A closer look at the above  proof reveals that
$Q(F^2(H_n)\otimes \cF)=F^2(H_n)\otimes \cF'$ and $V$ is a unitary operator. Taking into account  relation
\eqref{XX'} and \eqref{Th2}, we obtain
\begin{equation*}
\begin{split}
\Theta_2f\oplus X_\Theta^*(\Delta_2f\oplus 0)&=
\left[ \Theta_2'f'\oplus X_\Theta'^*(\Delta_2' f'\oplus 0)\right]+
\left[ 0\oplus X_\Theta'^*(0\oplus g')\right]\\
&=
\left[ \Theta_2'f'\oplus X_\Theta'^*(\Delta_2' f'\oplus 0)\right]+
\left[ 0\oplus X_\Theta^*(0\oplus V^*g')\right].
\end{split}
\end{equation*}
Hence,  we get
$$
\Theta_2 f\oplus X_\Theta^*(\Delta_2 f\oplus (-V^*g'))=
\Theta_2'f'\oplus X_\Theta'^*(\Delta_2' f'\oplus 0).
$$
Taking the norms, we have 
$$
\|f\|^2+\|g'\|^2=\|f'\|^2.
$$
Combining this with $\|f\|^2=\|f'\|^2+\|g'\|^2$, we obtain
$\|f\|=\|f'\|$, which shows that $Q$  is a unitary multi-analytic operator.
Due to \cite{Po-analytic}, this implies $Q=I\otimes \Psi_0$, for some unitary operator
$\Psi_0\in B(\cF,\cF')$.
Using relation
\eqref{first}, we complete the proof.
\end{proof}

\bigskip

\section{Triangulations for row contractions and joint invariant subspaces}

In this section, we prove the existence of a unique triangulation of type
\begin{equation}\label{C0-C1}
\left(\begin{matrix}C_{\cdot 0}&0\\
*& C_{\cdot 1}\end{matrix}
\right)
\end{equation}
for any row contraction $T:=[T_1,\ldots, T_n]$,   and prove the existence of joint invariant subspaces
for certain classes
of  row contractions.

We need a few  definitions.
Let $T:=[T_1,\ldots, T_n]$, $T_i\in B(\cH)$, be a row contraction. We say that  $T$ is of class $C_{\cdot 0}$ (or pure row contraction) if
$$
\lim_{k\to\infty}\sum_{|\alpha|=k}\|T_\alpha^* h\|^2=0\quad \text{for any }\quad h\in \cH,
$$
and of class
$C_{\cdot 1}$ if
$$
\lim_{k\to\infty}\sum_{|\alpha|=k}\|T_\alpha^* h\|^2\neq 0\quad \text{for any }\quad h\in \cH, ~h\neq 0.
$$
 We say that a row contraction $T:=[T_1,\ldots, T_n]$, $T_i\in B(\cH)$, has a triangulation of type \eqref{C0-C1}
if there is an orthogonal   decomposition $\cH=\cH_0\oplus \cH_1$ with respect to which
$$
T_i=\left(\begin{matrix} A_i&0\\
*& B_i\end{matrix}
\right),\quad i=1,\ldots,n,
$$
and the entries have the following properties:
\begin{enumerate}
\item[(i)] $T_i^*\cH_0\subset \cH_0$ for any $i=1,\ldots,n$;
\item[(ii)] $A:=[A_1,\ldots, A_n]$ is of class
$C_{\cdot 0}$;
\item[(iii)]
$B:=[B_1,\ldots, B_n]$ is of class
$C_{\cdot 1}$.
\end{enumerate}
The type of the entry denoted by $*$ is not specified.

\begin{theorem}\label{factori}
Every row contraction $T:=[T_1,\ldots, T_n]$, $T_i\in B(\cH)$, has a triangulation of type 
$$
\left(\begin{matrix}C_{\cdot 0}&0\\
*& C_{\cdot 1}\end{matrix}
\right)
$$
Moreover, this triangulation is uniquely determined.
\end{theorem}
\begin{proof}
First, notice that the subspace 
$$
\cH_0:=\left\{ h\in\cH:\ \lim_{k\to\infty} \sum_{|\alpha|=k}
\|T_\alpha^* h\|^2=0\right\}
$$
is invariant under each operator $T_i^*$, $i=1,\ldots,n$.
The decomposition
$\cH=\cH_0\oplus \cH_1$, where $\cH_1:=\cH\ominus \cH_0$, yields the triangulation
$$
T_i^*=\left(\begin{matrix} A_i^*&*\\
0& B_i^*\end{matrix}
\right),\quad i=1,\ldots,n,
$$
where $A_i^*:=T_i^*|_{\cH_0}$ and $B_i^*:=P_{\cH_1}T_i^*|_{\cH_1}$ for each $i=1,\ldots,n$.
Since
$$
\lim_{k\to\infty} \sum_{|\alpha|=k}
\|A_\alpha^* h\|^2=\lim_{k\to\infty} \sum_{|\alpha|=k}
\|T_\alpha^* h\|^2=0,\quad h\in \cH_0,
$$
the row contraction $A:=[A_1,\ldots, A_n]$ is of class
$C_{\cdot 0}$.
Now,
we need to show that 
$$
\lim_{k\to\infty} \sum_{|\alpha|=k}
\|B_\alpha^* h\|^2\neq 0 \quad \text{ for all }\quad h\in\cH_1, h\neq 0.
$$
Lt $V:=[V_1,\ldots, V_n]$, $V_i\in B(\cK)$, be the minimal isometric dilation of the row contraction $T:=[T_1,\ldots, T_n]$ (see Section 2).
For every $m=1,\ldots$, the isometries $V_\alpha$, $|\alpha|=m$, have orthogonal ranges.
Therefore, we have 
\begin{equation*}
\begin{split}
\left\|\sum_{|\alpha|=m}  V_\alpha \left(\sum_{|\beta|=k}
V_\beta T_\beta^* \right)P_{\cH_0}T_\alpha^* h
\right\|^2&=
\sum_{|\alpha|=m}\left\| \left(\sum_{|\beta|=k}
V_\beta T_\beta^* \right)P_{\cH_0}T_\alpha^* h
\right\|^2\\
&=
\sum_{|\alpha|=m} \sum_{|\beta|=k}
\left\| T_\beta^* P_{\cH_0}T_\alpha^* h
\right\|^2
\end{split}
\end{equation*}
for any $h\in \cH$.
  Since $P_{\cH_0}T_\alpha^* h\in \cH_0$,  we have
  \begin{equation}\label{Tb}
\lim\limits_{k\to\infty}
 \sum_{|\beta|=k}
\left\| T_\beta^* P_{\cH_0}T_\alpha^* h
\right\|^2=0.
\end{equation}
According to \cite{Po-isometric}, we have
\begin{equation}\label{PR1}
P_\cR h=\lim_{k\to\infty} \sum_{|\alpha|=k}
V_\alpha T_\alpha^* h,\quad \text{ for any } \quad h\in \cH,
\end{equation}
where $P_\cR$ is the orthogonal projection of the minimal isometric dilation space $\cK$ on the subspace $\cR$ in the Wold decomposition
$\cK=\cR\oplus M_V(\cL_*)$.
Now, using relations \eqref{Tb} and \eqref{PR1}, we obtain
\begin{equation*}
\begin{split}
P_\cR h
&=\lim_{k\to\infty} \sum_{|\alpha|=m}\sum_{|\beta|=k}
V_\alpha V_\beta T_\beta^* T_\alpha^* h\\
&=
\lim_{k\to\infty} \sum_{|\alpha|=m}  V_\alpha \left(\sum_{|\beta|=k}
V_\beta T_\beta^* \right)P_{\cH_0}T_\alpha^* h
+\lim_{k\to\infty} \sum_{|\alpha|=m}  V_\alpha \left(\sum_{|\beta|=k}
V_\beta T_\beta^* \right)P_{\cH_1}T_\alpha^* h\\
&=
  \sum_{|\alpha|=m}  V_\alpha  P_\cR P_{\cH_1}T_\alpha^* h.
\end{split}
\end{equation*}
Hence, we deduce that
\begin{equation*}
\begin{split}
\|P_\cR h\|^2&=\left\| \sum_{|\alpha|=m}  V_\alpha  P_\cR P_{\cH_1}T_\alpha^* h\right\|^2
=\sum_{|\alpha|=m}  \|  P_\cR P_{\cH_1}T_\alpha^* h\|^2\\
&\leq
\sum_{|\alpha|=m}  \|  P_{\cH_1}T_\alpha^* h\|^2
=\sum_{|\alpha|=m}  \|  B_\alpha^* h\|^2
\end{split}
\end{equation*}
for any $h\in \cH$.
Let $h\in \cH_1$, $h\neq 0$, and assume that 
$\lim\limits_{m\to\infty}\sum\limits_{|\alpha|=m}  \|  B_\alpha^* h\|^2=0$.
The above relation shows that $P_\cR h=0$ and, due to \eqref{PR1}, we deduce that
$h\in \cH_0$, which is a contradiction.

Now, we prove the uniqueness. Assume that
there is another  decomposition
$\cH=\cH_0'\oplus \cH_1'$ which yields the triangulation
$$
T_i=\left(\begin{matrix} {A_i'}&0\\
*& {B_i'}\end{matrix}
\right),\quad i=1,\ldots,n,
$$
  of type 
$
\left(\begin{matrix}C_{\cdot 0}&0\\
*& C_{\cdot 1}\end{matrix}
\right),
$
where ${A_i'}^*:=T_i^*|_{\cH_0'}$ and ${B_i'}^*:=P_{\cH_1'}T_i^*|_{\cH_1'}$ for each $i=1,\ldots,n$.
To prove uniqueness, it is enough to show that $\cH_0=\cH_0'$.
Notice that if $h\in \cH_0'$, then, due to the fact that
the row contraction
$[A_1',\ldots, A_n']$ is of class $C_{\cdot 0}$, we have
$$
\lim_{m\to\infty}\sum_{|\alpha|=m}  \|  T_\alpha^* h\|^2=\lim_{m\to\infty}\sum_{|\alpha|=m}  \|  {A_\alpha'}^* h\|^2=0.
$$
Hence, $h\in \cH_0$, which proves that $\cH_0'\subseteq \cH_0$.
Assume now that $h\in\cH_0\ominus \cH_0'$.
Since $h\in\cH_1'$, we have
$$
\lim\limits_{m\to\infty}\sum_{|\alpha|=m}  \|  {B_\alpha'}^* h\|^2=\lim\limits_{m\to\infty}\sum_{|\alpha|=m}  \|  P_{\cH_1'}{T_\alpha}^* h\|^2\leq
\lim\limits_{m\to\infty}\sum_{|\alpha|=m}\|  T_\alpha^* h\|^2=0.
$$
 Consequently,  since  the row contraction $[B_1',\ldots, B_n']$ is of class $C_{\cdot 1}$, we must have $h=0$. Hence, we deduce that $\cH_0\ominus\cH_0'=\{0\}$, which shows that $\cH_0'=\cH_0$. This completes the proof.
\end{proof}

\begin{corollary} If  $T:=[T_1,\ldots, T_n]$ is a row contraction such $T\notin C_{\cdot 0}$ and $T\notin C_{\cdot 1}$, then there is a non-trivial joint invariant subspace
under $T_1,\ldots, T_n$.
\end{corollary}

According to Section 2, any row contraction admits a triangulation of type
$$
\left(\begin{matrix}C_{c}&0\\
*& C_{cnc}\end{matrix}
\right)
$$
where $C_c$ (resp. $C_{cnc}$) denotes the class of coisometric (resp. c.n.c.) row contractions. Notice that $C_c\subset C_{\cdot 1}$.
Combining this result with the triangulation of Theorem
\ref{factori}, we obtain another triangulation for row contractions, that is,
$$
\left(
\begin{matrix}
C_{\cdot 0} &0&0\\
*& C_c &0\\
*&*& C_{cnc} \cap C_{\cdot 1}
\end{matrix}
\right).
$$

\begin{corollary}
If $T:=[T_1,\ldots, T_n]$, $T_i\in B(\cH)$,  is a row contraction such
$$T_1T_1^*+\cdots +T_nT_n^*\neq I$$
and there is a non-zero  vector $h\in \cH$ such that 
$\sum\limits_{|\alpha|=k}  \|  T_\alpha^* h\|^2=\|h\|^2$ for any $k=1,2,\ldots$, then there is a non-trivial subspace under the operators $T_1,\ldots, T_n$.
\end{corollary}

We recall from \cite{Po-similarity} that if 
$$T_1T_1^*+\cdots +T_nT_n^*=I,$$
then a subspace $\cM$ is invariant under $T_1,\ldots, T_n$ if and only if 
$$
T_1 P_\cM T_1^*+\cdots +T_nP_\cM T_n^*\leq P_\cM,
$$
where $P_\cM$ is the orthogonal projection on $\cM$.
We also  mention that
the case when $T\in C_{\cdot 0}$  is treated in the next corollary, and the case $T\in C_{\cdot 1}$ is considered in the next section
(see Theorem \ref{quasi-si}).

The proof of  the following result on regular factorizations of multi-analytic operators is straightforward from the definition, so we leave it to the reader.

\begin{lemma}\label{fact-reg}
Let $\Theta:F^2(H_n)\otimes\cE\to F^2(H_n)\otimes\cE_*$ be a 
contractive multi-analytic operator and assume that it has the factorization
$$
\Theta=\Theta_2\Theta_1,
$$
where $\Theta_1:F^2(H_n)\otimes\cE\to F^2(H_n)\otimes\cF$ and 
$\Theta_2:F^2(H_n)\otimes\cF\to F^2(H_n)\otimes\cE_*$
are contractive multi-analytic operators.
\begin{enumerate}
\item[(i)] If $\Theta_2$ is inner, then the factorization
$\Theta=\Theta_2\Theta_1$ is regular.
\item[(ii)]
If $\Theta$ is inner, then the factorization
$\Theta=\Theta_2\Theta_1$ is regular if and only if $\Theta_1$ and $\Theta_2$ are inner multi-analytic operators.
\item[(iii)]
If $\rank \Delta_\Theta<\infty$, then
$$
\rank \Delta_\Theta=\rank \Delta_{\Theta_2}+\rank \Delta_{\Theta_1}
$$
if and only if the factorization
$\Theta=\Theta_2\Theta_1$ is regular.
\end{enumerate}
\end{lemma}

Now we consider the case when $T$ is a pure row contraction.

\begin{corollary} \label{inner-inv}
If $T:=[T_1,\ldots,T_n]$ is a row contraction of class $C_{\cdot 0}$, then the non-trivial joint invariant subspaces under $T_1,\ldots, T_n$ are parametrized by the non-trivial inner factorizations of the characteristic function  $\Theta_T$ of T (i.e., $\Theta_T= \Theta_2\Theta_1$ with $\Theta_1$ and  $\Theta_2$ inner multi-analytic operators). Moreover, the subspaces $\HH_1$ and $\HH_2$ in Theorem $\ref{inv-factor1}$ become
\begin{equation*}
\begin{split}
\HH_1&=\{\Theta_2 f:\ f\in F^2(H_n)\otimes \cF\}\ominus \{\Theta_T f:\ f\in F^2(H_n)\otimes \cD\} \text{ and }\\
\HH_2&=\{F^2(H_n)\otimes \cD_*\}\ominus\{\Theta_2 f:\ f\in F^2(H_n)\otimes \cF\},
\end{split}
\end{equation*}
where $\cD$ and $\cD_*$ are the defect spaces of $T$.
\end{corollary}

\begin{proof}
According to Theorem \ref{funct-model1}, the characteristic function $\Theta_T$ is an inner multi-analytic operator. By Lemma \ref{fact-reg}, any factorization $\Theta_T=\Theta_2 \Theta_1$ is regular if and only if $\Theta_1$ and $\Theta_2$ are inner  operators. Applying now Theorem \ref{inv-factor1}, in our particular case, the result follows.
\end{proof}

We should remark that Corollary \ref{inner-inv} can also be proved directly using Theorem \ref{funct-model1} and the Beurling type characterization  (see \cite{Po-charact})
 of the joint invariant subspaces under   the operators $S_1\otimes I_\cG,\ldots, S_n\otimes I_\cG$.

We recall \cite{Po-multi} that any multi-analytic operator 
admits an essentially unique inner-outer factorization.

\begin{theorem}\label{chara-fact}
Let $T:=[T_1,\ldots, T_n]$ be a completely non-coisometric row contraction. The inner-outer factorization of the characteristic function $\Theta_T$ induces (cf.~Theorem $\ref{inv-factor2}$) the triangulation of type
$$
\left(\begin{matrix}C_{\cdot 0}&0\\
*& C_{\cdot 1}\end{matrix}
\right)
$$
for the row contraction $T$.

In particular,  if the inner-outer factorization of the characteristic function is non-trivial, then there is a non-trivial joint invariant subspace under the operators $T_1,\ldots, T_n$.
\end{theorem}

\begin{proof}
Suppose that the multi-analytic operator $\Theta:F^2(H_n)\otimes \cE\to F^2(H_n)\otimes \cE_*$ coincides with the characteristic function 
of the c.n.c. row contraction $T:=[T_1,\ldots, T_n]$.
Let $\Theta=\Theta_i\Theta_o$ be the cannonical inner-outer 
factorization of $\Theta$.
Since $\Theta_i$ is inner, Lemma \ref{fact-reg} implies that
the factorization is regular. Therefore, according to Theorem \ref{inv-factor1} (see also Theorem \ref{inv-factor}) and Theorem \ref{inv-factor2}, the above factorization yields a triangulation
$$
{\bf T}_i=\left(\begin{matrix}{\bf B}_i&0\\
*&{\bf A}_i
\end{matrix}\right),\quad i=1,\ldots,n,
$$
of ${\bf T}:=[{\bf T}_1,\ldots, {\bf T}_n]$, the functional model of $T$, such that the characteristic functions of  ${\bf B}:=[{\bf B}_1,\ldots, {\bf B}_n]$ and 
${\bf A}:=[{\bf A}_1,\ldots, {\bf A}_n]$ coincide with the purely contractive parts of $\Theta_i$ and $\Theta_o$, respectively. Due to Lemma \ref{pure}, the purely contractive part of an outer or inner multi-analytic operator is also outer or inner, respectively.
We recall from \cite{Po-charact} that a c.n.c. row contraction is of class $C_{\cdot 0}$ (resp. $C_{\cdot 1}$) if and only if the corresponding  characteristic function is inner (resp. outer) multi-analytic operator. Finally, using the last part of Theorem \ref{inv-factor2}, we  can complete the proof.
\end{proof}

\bigskip

\section{Characteristic functions and joint similarity to Cuntz row isometries}

In this section, we obtain criterions for joint similarity of $n$-tuples of operators to Cuntz row isometries. In particular, we prove that a  completely non-coisometric row contraction $T:=[T_1,\ldots, T_n]$ is jointly similar to a Cuntz row isometry if and only if the characteristic function of $T$  is an invertible multi-analytic operator.
This is a multivariable version of a result of Sz.-Nagy and Foia\c s \cite{SzF3}, concerning the  similarity to unitary operators.

 Extending on some results
obtained  by Sz.-Nagy   \cite {SzN}, Nagy-Foia\c s  \cite{SzF-book},  and the author \cite{Po-models}, \cite{Po-similarity}, we  provide necessary and sufficient conditions for a power bounded $n$-tuple  
 of operators  on a Hilbert space  to be  jointly similar to a Cuntz row isometry.

We need the following well-known result (see eg. \cite{SzF-book}).

\begin{lemma}\label{projec}
Let $\cM, \cN, \cX$ and $\cY$ be subspaces of a Hilbert space $\cH$ such that
$$
\cH=\cM\oplus \cN=\cX\oplus \cY.
$$
If 
$$P_\cM\cX=\cM \quad \text{ and } \quad \|P_\cM x\|\geq c\|x\|, \quad x\in \cX,
$$
 for some constant $c>0$, then 
$$
P_\cN\cY=\cN \quad \text{ and } \quad \|P_\cN y\|\geq c\|y\|, \quad y\in \cY.
$$
\end{lemma}

We recall  a few facts concerning the geometric structure of the minimal isometric dilation of a row contraction.
Let $T:=[T_1,\ldots, T_n]$, $T_i\in B(\cH)$, be a   row contraction
and let $V:=[V_1,\ldots, V_n]$ be its minimal isometric dilation   on a Hilbert space $\cK\supseteq \cH$. In \cite{Po-isometric}, we proved  that
$\cK=\cR\oplus M_V(\cL_*)$ and
\begin{equation}
\label{PR}
P_\cR h=\lim_{k\to\infty}\sum_{|\alpha|=k}V_\alpha T_\alpha^*h,\quad h\in \cH,
\end{equation}
where $P_\cR$ is the orthogonal projection of $\cK$ onto $\cR$.
Moreover, if $T$ is a one-to-one row contraction, then
\begin{equation}
\label{one-one}
\overline{P_\cR \cH}=\cR.
\end{equation}

The next result provides necessary and sufficient conditions for a c.n.c. row contraction to be jointly similar to a Cuntz row 
isometry, in terms of the corresponding characteristic function.

\begin{theorem}\label{simi-Cuntz}
Let $T:=[T_1,\ldots, T_n]$, $T_i\in B(\cH)$, be a completely
non-coisometric row contraction.  Then $T$ is jointly
similar  to a Cuntz row isometry
$W:=[W_1,\ldots, W_n]$, $W_i\in B(\cW)$, i.e.,
\begin{enumerate}
\item[(i)]
$W_1 W_1^*+\cdots +W_nW_n^*=I_\cW;
$
\item[(ii)] $ST_i =W_iS$, $i=1,\ldots, n$, for some invertible  operator $S:\cH\to \cW$,
\end{enumerate}
if and only if  the characteristic function $\Theta_T$
 is an invertible multi-analytic operator.

In this case, 
$$
\|\Theta_T^{-1}\|=\min \left\{\|X\|\|X^{-1}\|:\ [X^{-1} T_1 X,\ldots, X^{-1} T_n X] \ \text{ is a Cuntz row isometry}\right\}.
$$
\end{theorem}

\begin{proof}
Suppose that the row contraction
$T:=[T_1,\ldots, T_n]$ is  jointly similar  to a Cuntz row isometry
$W:=[W_1,\ldots, W_n]$, $W_i\in B(\cW)$, i.e.,
$$
W_1 W_1^*+\cdots +W_nW_n^*=I_\cW
$$
and  $T_i =S^{-1}W_iS$, $i=1,\ldots, n$, for some invertible  operator $S:\cH\to \cW$.
Since $ST_\alpha=W_\alpha S$ and $T_\alpha^* S^*=S^*W_\alpha^*$ for any $\alpha\in \FF_n^+$, we have
\begin{equation*}
\begin{split}
S\left(\sum_{|\alpha|=k}T_\alpha T_\alpha^*\right) S^*
&=
\sum_{|\alpha|=k}W_\alpha SS^*W_\alpha^*\\
&\geq \frac{1}{\|{S^*}^{-1} S^{-1}\|}\sum_{|\alpha|=k}W_\alpha W_\alpha^*\\
&=\frac{1}{\|S^{-1}\|^2} I
\end{split}
\end{equation*}
for any $k=1,2,\ldots$.
Therefore,
\begin{equation*}
\begin{split}
\sum_{|\alpha|=k}\left<T_\alpha T_\alpha^* h,h\right>
&\geq \|{S^*}^{-1}h\|^2\frac{1}{\|S^{-1}\|^2}\\
&\geq \frac{1}{\|S^*\|^2\|S^{-1}\|^2} \|h\|^2,
\end{split}
\end{equation*}
which, due to relation \eqref{PR}, implies
\begin{equation}\label{PR-ine}
\|P_\cR h\|\geq \frac{1}{\|S\|\|S^{-1}\|} \|h\|, \quad h\in \cH.
\end{equation}
Notice that the operator  $[T_1, \ldots, T_n]~$ is one-to-one. Indeed, the relation
$$
S^{-1}W_1 S h_1+\cdots +S^{-1}W_nS h_n=0,\qquad h_i\in \cH,~ i=1,\ldots,n,
$$
implies
$$
W_1S h_1+\cdots +W_nS h_n=0.
$$
Since $W_i$ are isometries with orthogonal ranges, we have
$$
W_i S h_i=0, \qquad i=1,\ldots, n,
$$
whence $~h_i=0$, $i=1,\ldots,n$. Therefore $~[T_1,\ldots, T_n]~$ is one-to-one.
According to  \eqref{one-one}, we  have $~\overline{P_\cR\cH}=\cR$.
Due to relation \eqref{PR-ine}, the subspace $P_\cR\cH$ is closed.
Therefore,  $P_\cR\cH=\cR$ and the operator
$$X:=P_\cR|_\cH:\cH\to \cR$$
is invertible.
According to \eqref{PR}, we have
\begin{equation*}\begin{split}
V_i^*P_\cR h&=\lim_{k\to\infty} \sum_{|\alpha|=k} V_i^*V_\alpha T_\alpha^*h\\
&=
\lim_{k\to\infty} \sum_{|\alpha|=k-1} V_\beta T_\beta^* T_i^*h=P_\cR T_i^*h
\end{split}
\end{equation*}
for any $h\in \cH$ and $i=1,\ldots, n$. 
 Consequently, we have
$$
 T_i X^*=X^*W_i,\qquad i=1,\ldots, n,
$$
where $W_i:=V_i|_\cR$,  $i=1,\ldots,n$.
Due to  the  noncommutative Wold decomposition applied to the row isometry
$[V_1,\ldots, V_n]$,  the subspace  $\cR$ is  reducing
  under each isometry $V_i$, $i=1,\ldots$, and 
$[W_1,\ldots, W_n]$  is a Cuntz row isometry.

Now, due to the geometric structure of the minimal isometric dilation of $T$ , we have (see  relation \eqref{two-dec})
$$
\cK=\cR\oplus M_V(\cL_*)=\cH\oplus M_V(\cL).
$$
Since $P_\cR \cH=\cR$, we can use   relation 
 \eqref{PR-ine} and Lemma \ref{projec} to deduce that
$$
P_{M_V(\cL_*)}M_V(\cL) =M_V(\cL_*)\quad \text{ and }\quad
\|P_{M_V(\cL_*)} x\|\geq \frac {1} {\|S\|\|S^{-1}\|}\|x\|,\quad x\in M_V(\cL).
$$
Therefore, the operator
$$
Q:=P_{M_V(\cL_*)}|_{M_V(\cL)}:M_V(\cL)\to M_V(\cL_*)
$$
is an invertible contraction  with ~$\|Q^{-1}\|\leq \|S\|\|S^{-1}\|$.
Since $Q$ is unitarily equivalent to the characteristic function $\Theta_T$ of $T$ (see Section 2),  we deduce that $\Theta_T$ is an invertible  multi-analytic operator and $\|\Theta_T^{-1}\|\leq \|S\|\|S^{-1}\|$.

Conversely, assume that the characteristic function $\Theta_T$ (and  hence $Q$) is an invertible contraction
and $\|\Theta_T^{-1}\|\leq \frac {1} {c}$ for some constant $c>0$.
Applying again Lemma \ref{projec}, we deduce that
$$
P_\cR\cH=\cR\quad \text{ and }\quad \|P_\cR h\|\geq c\|h\|,\quad h\in \cH.
$$
This shows that the operator $X:=P_\cR|_\cH:\cH\to \cR$ is invertible and $\|X^{-1}\|\leq \frac {1} {c}$.
As in the first part of the proof, we have
$X^*(V_i|_\cR)=T_iX^*$ for any $i=1,\ldots,n$.
This proves the similarity to a Cuntz row isometry. 
Notice also that, since $\|X\|\leq 1$, we have
$$
\|{X^*}^{-1}\|\|X^*\|=\|X^{-1}\|\|X\|\leq \frac {1} {c}.
$$

To prove the last part of the theorem,
let $c>0$ be such that $\|\Theta_T^{-1}\|=\frac {1} {c}$.
The converse of this theorem implies the existence of on invertible operator $X$ such that
 $[X^{-1} T_1 X,\ldots, X^{-1} T_n X]$ is a Cuntz row isometry and  
$$\|X\|\|X^{-1}\|\leq \frac {1} {c}=\|\Theta_T^{-1}\|.
$$
On the other hand, using the first part of the proof, we have
$$\|\Theta_T^{-1}\|\leq \|X\|\|X^{-1}\|.
$$
Therefore, $\|\Theta_T^{-1}\|= \|X\|\|X^{-1}\|$ and 
the proof is complete.
\end{proof}

\begin{corollary}
If 
$T:=[T_1,\ldots, T_n]$, $T_i\in B(\cH)$, is  a completely
non-coisometric row contraction  jointly
similar  to a Cuntz row isometry, then  
  $T$ is jointly similar to the Cuntz part   in the Wold decomposition
of the minimal isometric dilation of $T$.
Moreover,  in this case,  $T$
  is similar to 
the model row contraction $C:=[C_1,\ldots, C_n]$, where for each $i=1,\ldots, n$, 
$$
C_i: \overline{\Delta_{\Theta_T}(F^2(H_n)\otimes \cD)}\to
\overline{\Delta_{\Theta_T}(F^2(H_n)\otimes \cD)}
$$
is defined by
$$
C_i(\Delta_{\Theta_T}f):=\Delta_{\Theta_T}(S_i\otimes I_{\cD})f,\quad f\in F^2(H_n)\otimes \cD,
$$
and $\Delta_{\Theta_T}:=\left(I-\Theta_T^*\Theta_T\right)^{1/2}$, where $\Theta_T$ is the characteristic function of $T$. 
\end{corollary}

\begin{proof} The first part of the theorem follows from the proof of Theorem \ref{simi-Cuntz}. Now, using the model theory for c.n.c row contractions (see Theorem \ref{funct-model1} and Theorem \ref{funct-model2}),
one can complete the proof.
\end{proof}

Now we consider the case when $T:=[T_1,\ldots, T_n]$ is an arbitrary row contraction.
\begin{theorem}
  Let $T:=[T_1,\ldots, T_n]$, $T_i\in B(\cH)$,  be a row contraction.  Then $T$ is  jointly similar  to a Cuntz row isometry
$W:=[W_1,\ldots, W_n]$, $W_i\in \cW$, 
if and only if 
$T$ is one-to-one and  the operator 
\begin{equation}\label{P-SOT}
P:=\left(\text{\rm SOT-}\lim_{k\to\infty}\sum\limits_{|\alpha|=k}T_\alpha T_\alpha^*\right)^{1/2}
\end{equation}
is invertible.

Moreover, if this is the case, then the row contraction
$T:=[T_1,\ldots, T_n]$ is jointly similar to   the Cuntz part $R:=[R_1,\ldots, R_n]$ in the Wold decomposition
of the minimal isometric dilation of $T$.
\end{theorem}
\begin{proof}  Assume $T$ is a similar to $W$, i.e., there exists an 
invertible operator  $~S:\cH\to\cW~$  such that $~T_i=S^{-1}W_i S$, \ $i=1,\ldots, n$. 
As in the proof of Theorem \ref{simi-Cuntz}, one can show  that the operator  $[T_1, \ldots, T_n]~$ is one-to-one.  
According to \eqref{one-one}, we have $~\overline{P_\cR\cH}=\cR$.
On the other hand,  due to relation \eqref{PR}, we deduce that
\begin{equation}\label{PR-norm}
\|P_\cR h\|^2=\lim_{k\to\infty}\sum\limits_{|\alpha|=k}
\|T_\alpha^*h\|^2=\|P h\|^2,\qquad h\in\cH,
\end{equation}
where  the operator $P$ is well-defined by \eqref{P-SOT}, due to the fact that $\left\{\sum\limits_{|\alpha|=k}T_\alpha T_\alpha^*\right\}_{k=1}^\infty$
is a decreasing  sequence  of positive operators.
Notice that, since $\{W_\alpha\}_{|\alpha|=k}$ are isometries with orthogonal ranges, we have
\begin{equation*}
\begin{split}
\sum\limits_{|\alpha|=k}\|T_\alpha^*h\|^2 &\ge \|S^{-1}\|^{-2}
\sum\limits_{|\alpha|=k}\|W^*_\alpha S^{*-1}h\|^2 \\
&=\|S^{-1}\|^{-2}\|S^{*-1}h\|^2 \ge (\|S^{-1}\|^{2}\|S\|^2)^{-1} \|h\|^2
\end{split}
\end{equation*}
for any $h\in \cH$.
Therefore
$$~\|P_\cR h\|^2=\|P h\|^2\ge(\|S^{-1}\|^2\|S\|^2)^{-1}\|h\|^2$$
for any $~h\in\cH$. Hence, it follows that the operators $~P~$ and $~P_\cR|_\cH~$ are
one-to-one and have closed ranges. Since $~\overline{P_\cR\cH}=\cR$, it is clear that the operator
$~X:\cH\to\cR~$ is invertible.

According to relation \eqref{PR}, we have
\begin{equation*}
V_i^*P_\cR h= \lim_{k\to\infty} \sum_{|\alpha|=k-1} V_\beta T_\beta^* T_i^*h=P_\cR T_i^*h
\end{equation*}
for any $h\in \cH$ and $i=1,\ldots, n$. 
 Consequently, we deduce that
\begin{equation}\label{XTRX}
X T_i^*=R_i^*X,\qquad i=1,\ldots, n,
\end{equation}
where $~X:=P_\cR|_\cH$ and $R_i:=V_i|_\cR$,\ $i=1,\ldots,n$.
 Therefore,  $T:=[T_1,\ldots, T_n]$ is jointly similar to $R:=[R_1,\ldots, R_n]$.

Conversely, assume that the row contraction  $~[T_1,\ldots, T_n]~$ is one-to-one and  the operator $~P~$ is invertible.  Then relation \eqref{PR-norm}
implies $~P_\cR|_\cH~$ is one-to-one and has closed range.
On the other hand, by \eqref{one-one}, we have 
$~\overline{P_\cR\cH}=\cR$.
Therefore, the operator $~X:=P_\cR|_\cH:\cH\to \cR$ is invertible and, 
due to  relation \eqref{XTRX}, the row contraction $[T_1,\ldots, T_n]$ is jointly similar to the 
Cuntz row isometry
$[V_1|_\cR,\ldots, V_n|_\cR]$.
The proof is complete. 
\end{proof}
We recall (\cite{Po-similarity}) that an $n$-tuple
$[T_1,\ldots, T_n]$,  of operators $T_i\in B(\cH)$, is power bounded
if there is a  constant $M>0$
such that
$$
\sum_{|\alpha|=k} \|T_\alpha^*h\|^2\leq M^2\|h\|^2,\quad h\in\cH,
$$
for any $k=1,2,\ldots$.

\begin{theorem}\label{quasi-si}
Let $[T_1,\ldots, T_n]$ be a one-to-one power bounded $n$-tuple of operators on a Hilbert space $\cH$ such that, for any non-zero element $h\in \cH$, $\sum\limits_{|\alpha|=k} \|T_\alpha^*h\|^2$ does not converges to $0$ as $k\to\infty$. Then there exists a Cuntz row isometry
$[W_1,\ldots, W_n]$, $W_i\in B(\cH)$,   such that
$$
T_iX=XW_i,\quad i=1,\ldots, n,
$$
for some one-to-one operator $X\in B(\cH)$  with range dense
in $\cH$.
\end{theorem}
\begin{proof}
For each $h\in \cH$, $h\neq 0$, denote
$$
c(h):=\inf_{k=1,2,\ldots}\left(\sum_{|\alpha|=k} \|T_\alpha^*h\|^2\right)^{1/2}.
$$
Since $[T_1,\ldots, T_n]$ is a power bounded $n$-tuple of operators, there is aconstant $M>0$ such that
\begin{equation}\label{PB}
\sum_{|\alpha|=k} \|T_\alpha^*h\|^2\leq M^2\|h\|^2,\quad h\in\cH,
\end{equation}
for any $k=1,2,\ldots$. If $c(h)=0$ and $\epsilon>0$, then there is $k_0$ such that
$$
\left(\sum_{|\alpha|=k_0} \|T_\alpha^*h\|^2\right)^{1/2}\leq 
\frac{\epsilon}{M}.
$$
Hence and using \eqref{PB}, we deduce that
 \begin{equation*}
\begin{split}
\sum_{|\alpha|=m+k_0} \|T_\alpha^*h\|^2
&=
\sum_{|\beta|=k_0}\left< T_\beta\left(\sum_{|\gamma|=m} T_\gamma T_\gamma^*\right) T_\beta^* h,h\right>\\
&\leq M^2\sum_{|\beta|=k_0}
\left< T_\beta  T_\beta^* h,h\right>\leq \epsilon^2
\end{split}
\end{equation*}
for any  $m\geq 0$.
Consequently, $\lim\limits_{k\to\infty}\sum\limits_{|\alpha|=k} \|T_\alpha^*h\|^2=0$,
which contradicts the hypothesis. Therefore, we must have
$c(h)\neq 0$ for any $h\in \cH$, $h\neq 0$.

Now, for each $h,h'\in \cH$, we define
$$
[h,h']:=\operatornamewithlimits{LIM}_{k\to\infty}\sum_{|\alpha|=k}\left< T_\alpha^*h, T_\alpha^* h'\right>,
$$
where $\operatornamewithlimits{LIM}$ is a Banach limit.
Due to the properties of the Banach limit, $[\cdot,\cdot]$
is a bilinear form on $\cH$ and we deduce that
$$
[h,h]:=\operatornamewithlimits{LIM}_{k\to\infty}\sum_{|\alpha|=k}\| T_\alpha^*h\|^2\geq c(h)^2>0 \quad \text{ if } h\in \cH, ~h\neq 0,
$$
and $[h,h]\leq M^2\|h\|^2$. Moreover,  we have
$$
[h,h]=\sum_{i=1}^n [T_i^*h,T_i^*h],\quad h\in \cH.
$$
Due to a well-known theorem on bounded hermitian forms,
there exists a self-adjoint operator $P\in B(\cH)$
such that
$$
[h,h']=\left<Ph,h'\right>\quad \text{ for any } h,h'\in\cH,
$$
and, due to the above considerations, we have 
\begin{equation}\label{>0}
0<\left<Ph,h\right><M^2\|h\|^2,\quad h\in \cH, ~h\neq 0.
\end{equation}
Now, we show that $P=\sum\limits_{i=1}^n T_iPT_i^*$.
Indeed, we have
\begin{equation*}
\begin{split}
\left<Ph,h\right>&=
\operatornamewithlimits{LIM}_{k\to\infty}\sum_{|\alpha|=k+1}\| T_\alpha^*h\|^2
=\operatornamewithlimits{LIM}_{k\to\infty}\sum_{i=1}^n\sum_{|\alpha|=k}\| T_\alpha^*T_i^*h\|^2\\
&=
\sum_{i=1}^n [T_i^*h,T_i^*h]=\sum_{i=1}^n 
\left<PT_i^*h,T_i^*h\right>\\
&=\sum_{i=1}^n \left<\sum_{i=1}^n T_iPT_i^*h,h\right>
\end{split}
\end{equation*}
for any $h\in \cH$, which proves our assertion.
Notice that  relation \eqref{>0} shows that
the operator
$X:=P^{1/2}$ is one-to-one and has range dense in $\cH$.
Since $\sum\limits_{i=1}^n \|XT_i^*h\|^2=\|Xh\|^2$ for any $h\in \cH$,  it is clear that
$$
\sum_{i=1}^n \|XT_i^*X^{-1}x\|^2=\|x\|^2
$$
for any $x$ in the domain on $X^{-1}$. Hence and due to the fact that
the domain on $X^{-1}$ is dense in $\cH$,  the operators
$V_i^*:=XT_i^*X^{-1}$, $i=1,\ldots,n$,  can be extended by continuity on $\cH$. Using the same notation for the corresponding extensions, we have
$$\sum_{i=1}^n \|V_i^*h\|^2=\|h\|^2,\quad h\in \cH,
$$
and 
$
V_i^*X=XT_i^*$,\ $ i=1,\ldots,n.
$
This shows that $[V_1,\ldots, V_n]$ is a co-isometry from $\cH^{(n)}$ to $\cH$ such that
$$
T_iX=XV_i, \quad i=1,\ldots,n.
$$
Assume now that $h_i\in \cH$ and $\sum\limits_{i=1}^n V_ih_i=0$.
Then $\sum\limits_{i=1}^n T_iXh_i=0$.
Since $[T_1,\ldots, T_n]$ and $X$ are one-to-one operators, we must have
$h_i=0$  for each $i=1,\ldots,n$.
Consequently, $[V_1,\ldots, V_n]$ is a  one-to-one 
co-isometry, and therefore a unitary operator from  
$\cH^{(n)}$ to $\cH$. This implies that
$V_1,\ldots, V_n$ are isometries on $\cH$ with
$V_1V_1^*+\cdots +V_nV_n^*=I_\cH.
$
The proof is complete.
\end{proof}

As a consequence of Theorem \ref{quasi-si}, we deduce the following criterion for joint similarity of  a   power bounded $n$-tuple of operators  to a Cuntz row isometry.

\begin{corollary} Let $[T_1,\ldots, T_n]$ be a one-to-one power bounded $n$-tuple of operators on a Hilbert space $\cH$. Then $[T_1,\ldots, T_n]$
is jointly similar to a Cuntz row isometry if and only if
there exists  a constant $c>0$
such that
\begin{equation} \label{simi3}
\sum\limits_{|\alpha|=k} \|T_\alpha^*h\|^2\geq c\|h\|^2,\quad h\in\cH, 
\end{equation}
for any $k=1,2,\ldots$.
 \end{corollary}

\begin{proof} The direct implication can be extracted from the proof of Theorem \ref{simi-Cuntz}.
Conversely, if condition \eqref{simi3} holds, then, using the proof of Theorem
\ref{quasi-si}, we have
$$
c(h)\geq \sqrt{c} \|h\|,\quad h\in \cH, ~h\neq 0.
$$
Moreover, the positive operator $P\in B(\cH)$ has the properties
$$
T_iP^{1/2}= P^{1/2} V_i, \quad i=1,\ldots, n,
$$
where $[V_1,\ldots, V_n]$ is a Cuntz isometry,
and 
$$
\left<Ph,h\right>\geq c\|h\|^2,\quad h\in \cH, ~h\neq 0.
$$
Since the latter inequality shows that $P^{1/2}$ is an invertible operator, the result follows.
\end{proof}

\enddocument